\documentclass[12pt]{amsart}

\usepackage{amsmath, amssymb, amsthm, amsfonts, amscd}
\usepackage{tikz, tikz-cd}

\makeatletter
\def\@tocline#1#2#3#4#5#6#7{\relax
  \ifnum #1>\c@tocdepth 
  \else
    \par \addpenalty\@secpenalty\addvspace{#2}%
    \begingroup \hyphenpenalty\@M
    \@ifempty{#4}{%
      \@tempdima\csname r@tocindent\number#1\endcsname\relax
    }{%
      \@tempdima#4\relax
    }%
    \parindent\z@ \leftskip#3\relax \advance\leftskip\@tempdima\relax
    \rightskip\@pnumwidth plus4em \parfillskip-\@pnumwidth
    #5\leavevmode\hskip-\@tempdima
      \ifcase #1
       \or\or \hskip 3em \or \hskip 6em \else \hskip 9em \fi%
      #6\nobreak\relax
    \hfill\hbox to\@pnumwidth{\@tocpagenum{#7}}\par
    \nobreak
    \endgroup
  \fi}
\makeatother

\input xy
\xyoption{all}

\usepackage{hyperref}
\usepackage{graphicx}
\usepackage{color}
\usepackage{verbatim}

\theoremstyle{plain}
\newtheorem{theorem}{Theorem}[subsection]
\newtheorem{lemma}[theorem]{Lemma}
\newtheorem{proposition}[theorem]{Proposition}
\newtheorem{prop}[theorem]{Proposition}
\newtheorem{cor}[theorem]{Corollary}

\theoremstyle{definition}
\newtheorem{defn}[theorem]{Definition}

\newtheorem{remark}[theorem]{Remark}

\newtheorem{example}[theorem]{Example}

\usepackage[margin=1in,marginparwidth=0.8in, marginparsep=0.1in]{geometry}

\def\AA{\mathbb{A}}

\def\CC{\mathbb{C}}

\def\G{\mathbb{G}}
\def\HH{\mathbb{H}}

\def\LL{\mathbb{L}}

\def\NN{\mathbb{N}}

\def\PP{\mathbb{P}}

\def\RR{\mathbb{R}}
\def\R{\mathbb{R}}

\def\TT{\mathbb{T}}

\def\ZZ{\mathbb{Z}}
\def\BZ{\mathbb{Z}}

\def\calA{\mathcal{A}}

\def\calC{\mathcal{C}}

\def\cF{\mathcal{F}}

\def\calP{\mathcal{P}}

\def\calT{\mathcal{T}}

\def\bT{\mathbf{T}}

\newcommand\tilcalP{\widetilde{\mathcal{P}}}

\def\hatMc{\widehat{M^\vee}}

\newcommand{\Coh}{\textup{Coh}}

\newcommand\Fuk{\textup{Fuk}}

\newcommand{\IndCoh}{\textup{IndCoh}}
\newcommand{\IndFuk}{\textup{IndFuk}}

\newcommand\Mod{\textup{Mod}}

\newcommand{\Perf}{\textup{Perf}}

\newcommand{\QCoh}{\textup{QCoh}}

\newcommand\Spec{\textup{Spec}}

\newcommand\Hom{\textup{Hom}}

\renewcommand{\L}{\Lambda}

\newcommand\nc{\newcommand}
\nc\on{\operatorname}
\nc\ol{\overline}
\nc\ul{\underline}

\nc\oo{\infty}

\nc\Cone{\mathit{Cone}}
\nc\ssupp{\mathit{ss}}
\nc\risom{\stackrel{\sim}{\to}}
\nc\Sh{\mathit{Sh}}
\nc\un{\diamondsuit}
\nc\orient{\mathit{or}}
\nc\sing{\mathit{sing}}
\nc\MF{\on{MF}}
\nc\Log{\on{Log}}
\nc\Arg{\on{Arg}}
\nc\inthom{\mathit{Hom}}
\nc\colim{\on{colim}}

\nc\wdv{W_{\Delta^\vee}}
\nc\trt{\widetilde{\on{trop}}}
\nc\dtmir{\partial\bT^{mir}}
\nc\ttmir{\partial\bT^{mir}_{loc}}
\nc\FS{\on{FS}}
\nc\tcP{\widetilde{\mathcal{P}}}
\nc\tcA{\widetilde{\mathcal{A}}}
\nc\Conv{\on{Conv}}
\nc\spb{\overline{\sigma}^\perp}
\nc\oS{\overline{S}}
\nc\Ain{\partial^0\tcA}

\nc\conv{\mathit{conv}}
\nc\trad{\mathit{inf}}
\nc\wrap{\mathit{wr}}

\newcommand{\ocfan}[1]{{\hat{\Sigma}(\overline{O(#1)})}}

\usepackage{tikz}
\usetikzlibrary{matrix,shapes,arrows,calc,topaths,intersections,positioning,decorations.pathreplacing,decorations.pathmorphing,fit}

\newcommand*{\hrlen}{10}
\newcommand*{\hramp}{3}

\tikzset{
righthairs/.style={postaction={decorate,draw,decoration={border,amplitude=\hramp,segment length=\hrlen,angle=-90,pre=moveto,pre length=\hrlen/2}}},
lefthairs/.style={postaction={decorate,draw,decoration={border,amplitude=\hramp,segment length=\hrlen,angle=90,pre=moveto,pre length=\hrlen/2}}},
righthairsnogap/.style={postaction={decorate,draw,decoration={border,amplitude=\hramp,segment length=\hrlen,angle=-90}}},
lefthairsnogap/.style={postaction={decorate,draw,decoration={border,amplitude=\hramp,segment length=\hrlen,angle=90}}},
}

\begin{document}


\title[Mirror symmetry for very affine hypersurfaces]{Mirror symmetry for very affine hypersurfaces}

\author{Benjamin Gammage and Vivek Shende}

\begin{abstract}
We show that the category of coherent sheaves on the toric boundary divisor of a smooth quasiprojective toric DM stack
is equivalent to the wrapped Fukaya category of a hypersurface in $(\CC^\times)^n$.  Hypersurfaces with every
Newton polytope can be obtained. 

Our proof has the following ingredients.  
Using recent results on localization, we may trade 
wrapped Fukaya categories for microlocal
sheaf theory along the skeleton of the hypersurface. 
Using Mikhalkin-Viro patchworking,
we identify the skeleton of the hypersurface with 
the boundary of the Fang-Liu-Treumann-Zaslow skeleton.
 By proving a new functoriality result for Bondal's coherent-constructible correspondence, we 
reduce the sheaf calculation to Kuwagaki's recent theorem on mirror symmetry for toric varieties.
\end{abstract}

\maketitle

\thispagestyle{empty}

\newpage 
\tableofcontents

\thispagestyle{empty}

\newpage

\section{Introduction}

Homological mirror symmetry is a story of two categories radically different in origin.  The first is a category of Lagrangians
in a symplectic manifold, with morphisms defined by intersection points, corrected by holomorphic disks.  The second
is a category of locally defined modules over the holomorphic functions on a seemingly unrelated complex variety, 
with morphisms corrected by considerations of homological algebra.  Most articles on the subject 
concern the ingenious manipulations required to identify one with the other, most often requiring heroic
calculations of at least one side of this equivalence. 

Our contribution is of a different nature.  We wish to explain how in many circumstances -- we focus on Calabi-Yau hypersurfaces in 
toric varieties, though the same methods should apply in the generality of Gross-Siebert toric degenerations -- both 
sides can be cut into matching elementary pieces, known to be homologically mirror, and the total mirror symmetry glued together using
foundational results in algebraic and symplectic geometry.  More precisely, this cutting and gluing is possible at the limiting point where on the one hand 
the complex manifold degenerates into a union of toric varieties, while on the other, the symplectic
form concentrates along certain divisors, and we consider the category associated to their complement.  We will
be entirely concerned with homological mirror symmetry {\em at} this limit point.\footnote{It is a tautology that matching the limit categories
 matches their infinitesimal deformations, but it remains to identify the geometric meaning of these 
deformations in a satisfactory way -- we do not touch upon this question here.} 

At this most degenerate point, the category of coherent sheaves on the union of toric varieties -- glued together 
along toric subvarieties -- can be calculated as a colimit
of the categories of coherent sheaves on the toric components \cite{GR}.  

Mirror symmetry is well studied for toric varieties themselves.  The Hori-Vafa prescription
is that the mirror A-model category should be associated to a function 
$W: (\CC^*)^n \to \CC$ whose Newton polytope is the 
convex hull of primitive vectors on the 1-dimensional cones of the fan of the toric variety.  
Different authors have taken different views
on how precisely to associate a category to this geometry, either directly in 
Lagrangian Floer theory \cite{A1, A2}, or in microlocal sheaf theory \cite{B, FLTZ, Tr, Ku}
(the latter being known to be calculate Fukaya categories \cite{NZ, Ncs, GPS3}). 

A true believer in mirror symmetry should expect the following facts: 

\begin{enumerate}
\item \label{boundarycovered} The mirror to the toric boundary -- a generic fiber of a 
generic $W: (\CC^*)^n \to \CC$ whose Newton polytope is the 
moment polytope of the toric variety -- admits a cover by mirrors of toric varieties, glued along 
mirrors of toric varieties. 
\item \label{boundaryfunctors} There are geometrically defined functors between these Fukaya categories which are mirror to the 
pullback and pushforward functors corresponding to the inclusion of toric varieties in toric varieties.
\item \label{descent} There is a descent result for the Fukaya category showing that it carries covers of the sort in 
(\ref{boundarycovered}) to colimits of categories. 
\end{enumerate} 

Establishing all of these results would show that the Fukaya category of the general fiber of $W$ 
is equivalent to the category of coherent sheaves of the corresponding toric variety.  
The recent works \cite{GPS1, GPS2, GPS3} give the necessary general tools to define the functors
in (\ref{boundaryfunctors}) and establish the descent required in (\ref{boundarycovered}).  In fact, these works, 
together with \cite{NS}, build a bridge between Fukaya categories and microlocal sheaf theory, which we 
cross in order to appeal to the microlocal sheaf calculations of the toric mirror \cite{Ku}.  
Here we will establish (\ref{boundarycovered}) and the `mirror' assertions of (\ref{boundaryfunctors}) 
above, and deduce:

\begin{theorem} \label{maintheorem} 
Suppose we are given the following data: 
\begin{itemize}
\item 
$\TT_\CC$  an algebraic torus with character and cocharacter lattices  
$M$ and $M^\vee$.  
\item  $\Delta^\vee \subset M^\vee$  an integral polytope containing
the origin.  
\item $\Sigma$  a fan in $M^\vee \otimes \RR$ giving a star-shaped triangulation
of $\Delta^\vee$. 
\end{itemize}
These determine a smooth toric stack  $\bT_\Sigma$ with toric boundary divisor $\partial \bT_\Sigma$. 

Then there exists a Laurent  polynomial $W: \TT^\vee_\CC \to \CC$ with Newton polytope $\Delta^\vee$;
a natural structure of Liouville manifold on a general fiber $F_W$;  and an equivalence
$$\Coh(\partial\bT_{\Sigma})\cong\Fuk(F_W)$$
	between the dg category of coherent sheaves on the variety $\partial\bT_\Sigma$ and the wrapped Fukaya category of the general fiber $F_W$
\end{theorem} 

We close the introduction with some comments about how we will establish items (\ref{boundarycovered})
and (\ref{boundaryfunctors}) above.  

Regarding 
(\ref{boundarycovered}): in the microlocal sheaf theoretic works beginning with \cite{FLTZ}, 
a key role is played by a certain conical Lagrangian subvariety 
$\Lambda_\Sigma \subset T^* \TT \cong \TT^\vee_\CC$.  It is straightforward to establish 
that the boundary of this conical subvariety indeed admits a cover corresponding to the cover of 
the $\partial \bT$ by toric subvarieties.  What is needed is to relate the geometry of $\Lambda_\Sigma$ 
to the geometry of the Laurent polynomial $W$.  One result
along these lines which we shall establish is that the deformation equivalence
class of the Liouville sector determined
by $W$ admits a representative whose relative skeleton is precisely $\Lambda_\Sigma$.  
Another 
is that a neighborhood of the boundary of this sector admits a sectorial cover by sectors whose relative 
skeleta give the aforementioned cover of $\Lambda_\Sigma$.  These results are established in
Section \ref{sec:skeleta} by using the 
Mikhalkin-Viro patchworking to reduce the study of $F_W$ to understanding pairs of pants, 
whose skeleta have been calculated by Nadler.

Regarding (\ref{boundaryfunctors}): after the geometric results in the previous paragraph, existence of
the relevant functors of Fukaya categories can be deduced from \cite{GPS1}.  To calculate them, 
we use \cite{GPS3, NS} to pass to microlocal sheaf theory, where we must now show that the
mirror symmetry established in \cite{Ku} can be made functorial with respect to inclusion of toric boundary 
divisors.  We explain how to do this in Section \ref{sec:microccc}.

In the following section, we explain in more detail the general strategy of the proof,
reviewing relevant ideas from sources mentioned above, and we give the proof of
Theorem \ref{maintheorem}, up to the calculations mentioned in the previous two paragraphs,
which we defer to the main body of the paper.

\bigskip
\emph{Acknowledgements} ---
We thank Sheel Ganatra, Stephane Guillermou, Allen Knutson, Tatsuki Kuwagaki, Heather Lee, Grigory Mikhalkin, David Nadler, Martin Olsson, John Pardon, Pierre Schapira, Nick Sheridan, Laura Starkston, and Zack Sylvan
for helpful conversations on various related topics. The authors are also thankful to Peng Zhou for notifying us of an omission in an earlier version of this paper. 

The work of B.G. was supported by an NSF Graduate Research Fellowship, and V.S. was supported by NSF DMS-1406871, NSF CAREER DMS-1654545, and a Sloan fellowship.

\section{Our approach to homological mirror symmetry}  \label{sec:discussion}

\subsection{An illustration} \label{pictures}

Consider the degeneration in which a genus-one holomorphic curve acquires a node.
In the mirror degeneration, a symplectic 2-torus acquires a puncture.

\begin{figure}[h]
\includegraphics[scale=0.5]{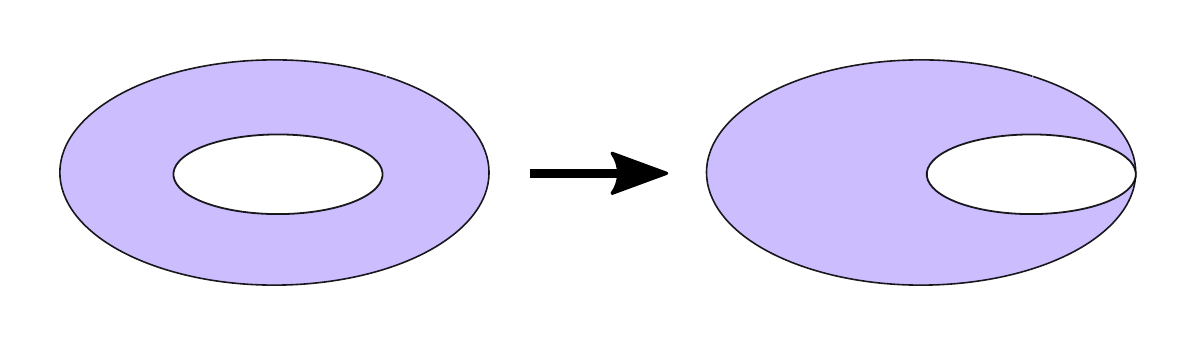}
\caption{The degeneration of a smooth genus-one curve to a nodal curve.}
\end{figure}

\begin{figure}[h]
\includegraphics[scale=0.5]{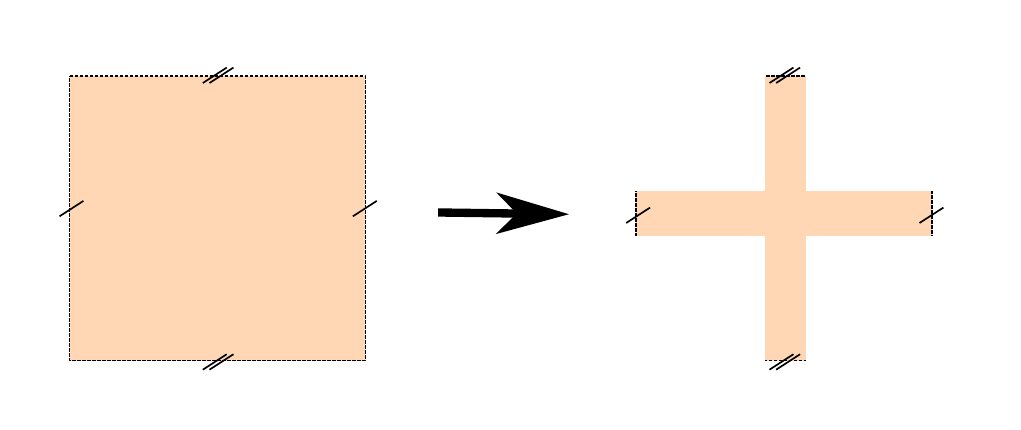}
\caption{\label{fig:puncture}A torus acquiring a puncture as an $S^1$ fiber approaches infinite radius.}
\end{figure}

\begin{figure}[h]
\includegraphics[scale=0.5]{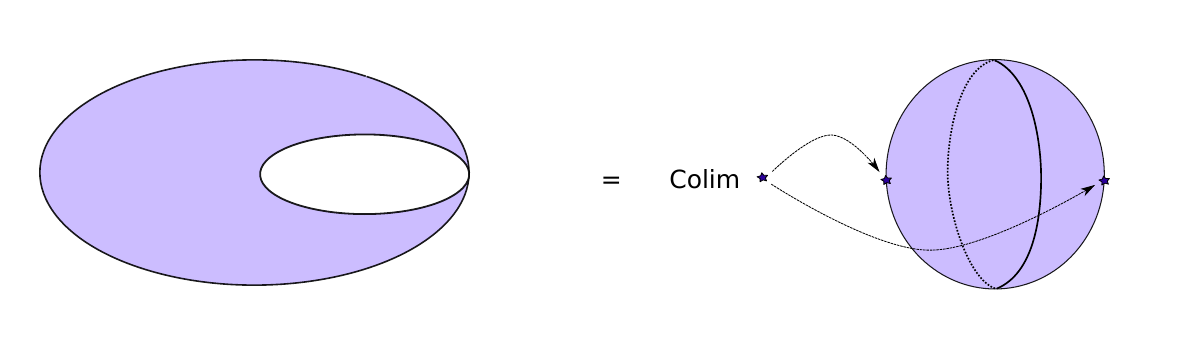}
\caption{\label{fig:normalization}We obtain a nodal curve by gluing smooth pieces.}
\end{figure}

\begin{figure}[h]
\includegraphics[scale=0.5]{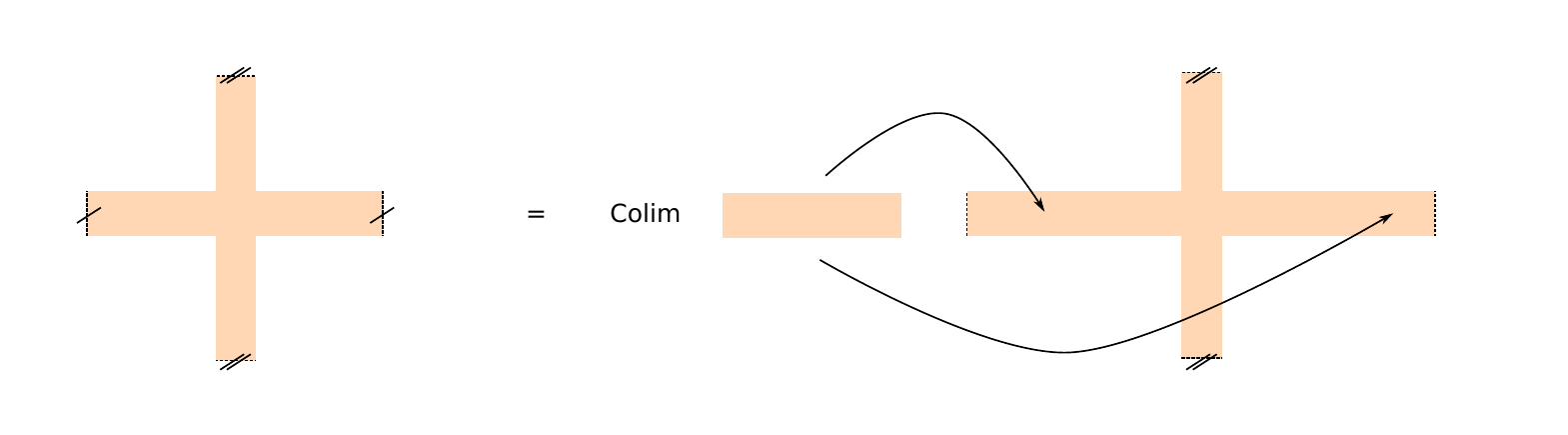}
\caption{ \label{fig:sectorgluing} The mirror to the above gluing: a punctured torus is glued together from Liouville sectors.}
\end{figure}

\newpage

One way to arrive at the view that these two spaces should be mirror is the following ``T-duality'' account.
In general, the spaces on the two sides of mirror symmetry are expected to be dual torus fibrations (in general, with singularities) over the same base, 
the radii of the fibers on one side being inverse to the radii on the other side.   In the present example, on the complex side,
we have a torus -- a circle bundle over a circle.  Under the degeneration, one of the circle fibers is approaching zero radius.
Thus, on the symplectic side, we should have a circle bundle over a circle, in which one fiber is approaching infinite radius.  A circle
of infinite radius is a line -- or in other words, the fiber should acquire a puncture.  

In the description above, the puncture was just the removal of a point.  As we draw only the complement of this point, we are free 
to imagine the puncture as being larger, as in Figure \ref{fig:puncture}.  
In our previous description, the fiber containing the puncture was dual to the node.  We have expanded the puncture, so in this
picture, one should regard the entire horizontal region beneath the puncture as being dual to the node. 

On the complex side, 
we have a singular complex curve; it is natural to take the  normalization.  This is a smooth curve mapping to the singular
curve, and in the case at hand, the map simply identifies points.  This is what is indicated in Figure \ref{fig:normalization}. 
We can describe the symplectic side by a similar gluing.  Since the node corresponded to the strip beneath
the puncture, the mirror gluing on the A-side involves gluing the two ends of the strip.  

The category we associate to this noncompact symplectic manifold is the wrapped Fukaya category, which was originally
constructed for \emph{Liouville manifolds}, symplectic manifolds with the property that (at least locally near the boundary) there is a primitive 
for the symplectic form whose dual Liouville vector field is everywhere outward pointing \cite{AS}.  
In the above gluing, however, the restriction of this
Liouville form to the components {\em does not} have this property: there are boundary components where it is
parallel, rather than outward pointing.  In particular, the rectangle should be viewed as the cotangent bundle of an interval
rather than a disk.  That is, the pieces in our gluing are not Liouville manifolds.   
The appropriate notion is that of {\em Liouville sector}, which we review in the next subsection. 
A covariantly functorial Floer theory for these is developed in \cite{GPS1, GPS2}.

\begin{figure}[h]
\includegraphics[scale=0.7]{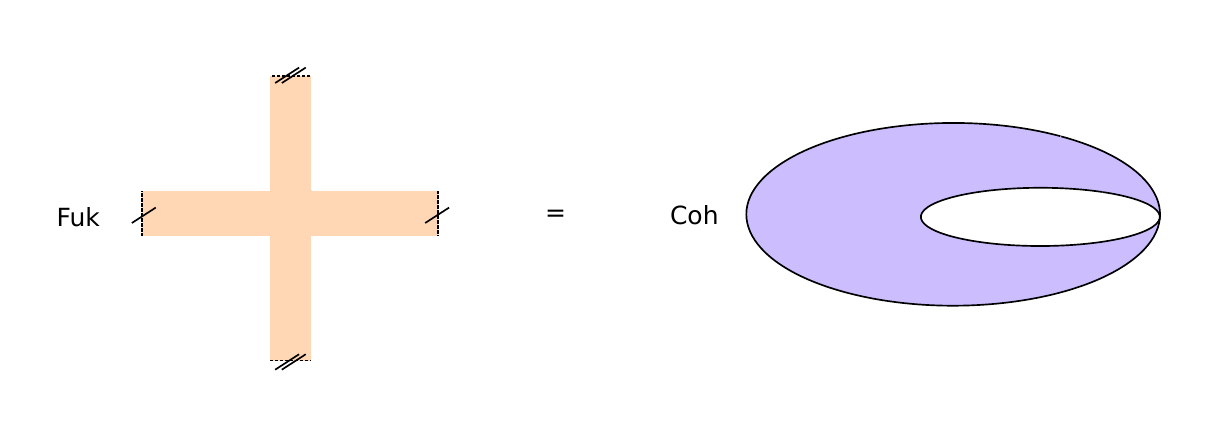}
\caption{\label{fig:mirrorexample}The homological mirror symmetry conjecture for a genus-one curve at the large volume/complex structure limits.}
\end{figure}

We turn now to the question of gluing together a global mirror symmetry from local mirror symmetries. 
The functor $\Coh(-)$ taking a variety $X$ to its dg category $\Coh(X)$ of coherent sheaves behaves well with respect to the gluings above.  
Following \cite{GR}, we make statements for the Ind-completion $\IndCoh(X)$ of the category $\Coh(X)$; 
statements for $\Coh(X)$ may be recovered by taking compact objects.  We write $\IndCoh^!$ and $\IndCoh_*$ for the contravariant and covariant
functors from derived stacks to dg categories which carry a stack to its
category of Ind-coherent sheaves and carry a map $f:X\to Y$ to a pullback $f^!$
or a pushforward $f_*$, respectively.  The key fact 
\cite[IV.4.A.1.2]{GR} is that
$\IndCoh^!$ takes pushout squares of affine\footnote{The above schemes are not affine, but the desired
pushout formula can be checked affine locally by Zariski descent.} 
schemes along closed embeddings to pullback squares of (stable cocomplete) dg categories.
By passing to adjoints, we see that $\IndCoh_*$ analogously takes pushouts to pushouts.

Homological mirror symmetry as usually stated is an equivalence between coherent sheaves on a given algebraic variety and 
the Fukaya category of its mirror.  What the pictures above suggest is that this should extend to a natural transformation
between the functor $\IndCoh_*$, perhaps with respect to some restricted class of maps including normalizations, and a
functor $\IndFuk_*$, covariant with respect to some class of maps including those mirror to normalization.  It suggests moreover that  $\IndFuk_*$ should take certain diagrams -- those mirror to certain pushouts of varieties -- to pushouts of dg or $A_\infty$ categories.

In fact, a covariant functor from a category of Liouville sectors to $A_\infty$ categories has been defined in \cite{GPS1}, and 
shown in \cite{GPS2} to carry diagrams like those illustrated above to pushouts.  
Given these structural properties, one can establish mirror symmetry by showing that 
there is an identification, respecting the relevant inclusion functors, of 
the Fukaya and coherent-sheaf categories of our building blocks.

\begin{remark}
	Strictly speaking, this subsection does not illustrate a special case of the statement of Theorem \ref{maintheorem}, 
	because the self-nodal curve is not the boundary of a toric variety.  However, an essentially 
	identical argument gives the mirror symmetry 
	between the nodal necklace of three $\mathbb{P}^1$s and a thrice-punctured torus, 
	which is a special case of Theorem \ref{maintheorem}. 
	
	This subsection also does not exactly illustrate the proof we will give of Theorem \ref{maintheorem}; instead, as we 
	explain below, we will translate these ideas to the microlocal sheaf setting using \cite{GPS3}.  
	In this setting, we need only cover the skeleton, as we do in Corollary \ref{cor:cover1}.  A lift of this cover to 
	a sectorial cover in the sense of \cite{GPS2} would yield a proof hewing closer to the above illustration. 
\end{remark}

\subsection{Stops, sectors, skeleta, and partially wrapped Fukaya categories.}
For basics on Liouville and Weinstein manifolds, we refer to \cite{CE, Eli-revisited}.  
Here we review basic notions of Liouville sectors, stops, and skeleta, and then recall 
from \cite{GPS1, GPS2} definitions and results concerning partially wrapped Fukaya categories defined
in terms of these geometric structures. 

A {\em Liouville sector} is an exact symplectic manifold-with-boundary $(X, \partial X, \lambda)$
modeled at infinity on the symplectization of a contact manifold-with-boundary $(V, \partial V, \lambda)$,
satisfying additional constraints: $\partial V$ should be transverse to a contact vector field, 
and the characteristic foliation on $\partial X$ should be trivializable as $\partial X = \mathbb{R} \times F$. 
Such a trivialization makes $(F, \lambda|_F)$ a Liouville manifold.  
Note that being a Liouville sector is a property of, rather than a structure on, an exact symplectic
manifold-with-boundary. 

A closed codimension zero 
submanifold-with-boundary $Y \subset X$ is a Liouville subsector if (1) each component of 
$\partial Y$ is either disjoint from or contained in $\partial X$, and (2) $(Y, \partial Y, \lambda|_{\partial Y})$
is itself a Liouville sector.   One can see by inspection that the symplectic manifolds in 
Figure \ref{fig:sectorgluing} admit an exact structure making them Liouville sectors, and that the inclusions
depicted are inclusions of Liouville sectors. 

Another point of view on sectors is obtained by passing to the `convex completion' $\bar{X}$, which is 
a Liouville manifold in the usual sense.  Up to contractible choices, the data of the sector is equivalent
to an embedding $F \subset \partial_\infty \bar{X}$ as a Liouville hypersurface, i.e. some choice of contact
form on $\partial_\infty \bar{X}$ restricts to the Liouville form on $F$.  The $(\bar{X}, F)$ form what
is termed a {\em Liouville pair} elsewhere in the literature \cite{Avdek, Sylvan, 
Eli-revisited}.  The advantage of Liouville sectors over Liouville pairs 
is that they are better suited to discussions of gluing, in particular because the key 
notion of Liouville subsector is less natural in the setting of pairs. 
Basic definitions and constructions relevant
to Liouville sectors are found in \cite[Sec. 2]{GPS1}.

We refer to works \cite{AS, GPS1, GPS2} for a foundational treatment of partially wrapped Fukaya
categories.  For our purposes here, we may largely use these works as black boxes.  The most
general setting for defining partially wrapped Fukaya categories (offered by \cite{GPS2}) takes
as input the data of a Liouville sector $X$ and a closed subset in the infinite boundary 
$\Lambda \subset \partial_\infty X^\circ$.  (We recall that a Liouville sector has its actual boundary 
$\partial X$, and its ideal contact boundary $\partial_\infty X$; here $\partial^\infty X^\circ$ means
the contact boundary minus its intersection with the actual boundary $\partial X$.)  To such a pair
is associated a category which we here denote $Fuk(X, \Lambda)$.\footnote{We write $Fuk$ for what is called $\mathrm{Perf}\, \mathcal{W}$ in \cite{GPS1, GPS2, GPS3}.}  

A stopped sector includes in another by enlarging the sector or shrinking the stop: we say 
$(X', \Lambda') \subset (X, \Lambda)$ if $X'$ is a Liouville subsector of $X$ and 
$\Lambda' \supset \Lambda \cap X'$.   It is shown in \cite{GPS1, GPS2} that in this case 
there is a functor
$$Fuk(X', \Lambda') \to Fuk(X, \Lambda)$$ 
When $X = X'$ we term
this functor a ``stop removal.'' 
These satisfy the natural compatibilities with composition, defining a (strict!) functor from 
the poset of (stopped) subsectors of $(X, \Lambda)$ to $A_\infty$ categories.

To a Liouville manifold $X$ is associated the {\em skeleton} (elsewhere termed {\em spine} or {\em core}) $\mathfrak{c}_X$, this being the locus 
of all points which do not escape to infinity under the Liouville flow.  When the Liouville flow
is gradient-like and generalized Morse-Smale (such manifolds are said to be Weinstein), the skeleton is 
admits a Whitney stratification by isotropic submanifolds, and the top-dimensional strata admit transverse ``cocore'' Lagrangian disks.  It is this consequence
which is relevant for \cite{GPS2, GPS3},\footnote{Without 
any assumption beyond isotropicity of the skeleton, one can use the linking
disks in $X \times T^*[0,1]$ as replacements for the co-core disks of $X$; see \cite{GPS2}.} 
and some weaker definitions of Weinstein have been proposed which imply it; see for instance \cite{Eli-revisited}.
The above results remain true when the Liouville flow is Morse-Bott, as in the cases studied in this paper.

For a Liouville sector $X$, one can define the skeleton $\mathfrak{c}_X$ 
by the same formulation: $\mathfrak{c}_X$ is the locus which does
not escape to infinity.  However, this definition is only really sensible if the Liouville flow on 
$X$ is tangent to $\partial X$ along all of $\partial X$, not just at the boundary.  
Note \cite[Lemma 2.11, Prop 2.28]{GPS1} this can always be arranged after
deformation.  Evidently, if 
$X \subset Y$ is an inclusion of Liouville sectors where $\lambda_X = \lambda_Y|_X$ 
and the Liouville flow on $X$ is tangent to its boundary, then $\mathfrak{c}_X = \mathfrak{c}_Y \cap X$.

We offer also another perspective on the skeleton of a sector.  Recall that a sector $X$ is equivalent to 
the data of a pair $(\bar{X}, F \subset \partial_\infty \bar{X})$.  For the pair, it is natural to define 
the {\em relative skeleton} $\mathfrak{c}_{\bar{X}, F}$ as the locus of points $\bar{X}$ 
which do not escape to $\partial_\infty \bar{X} \setminus \mathfrak{c}_F$.  Note that this
is the union of $\mathfrak{c}_X$ with a $\R$-cone on $\mathfrak{c}_F$.  This notion of relative skeleton
compares  to the skeleton of a sector as follows: it is not difficult, 
using the techniques of \cite[Sec. 2]{GPS1} to arrange an inclusion of sectors 
$X \subset \bar{X}$ such that $\mathfrak{c}_X = \mathfrak{c}_{\bar{X}, F} \cap X$. 

From the point of view of Fukaya categories, the significance of skeleta and relative skeleta is in their
role in organizing generation results.  Indeed, the cocore disks to a Weinstein Morse function 
provide Lagrangians transverse to each component of the smooth locus of the skeleton, and 
at any Legendrian point of a stop there is associated a linking disk; according 
to \cite[Thm. 1.10]{GPS2}, these generate $Fuk(X, \Lambda)$ when $X$ is a Weinstein manifold
and $\Lambda$ is mostly Legendrian.  

For calculating Fukaya categories, we may always translate back and forth between Liouville sectors
and stopped Liouville manifolds, and further we may retract the stop to its skeleton.  Indeed, 
per \cite[Cor. 2.11]{GPS2}, we have equivalences

$$Fuk(X) \xrightarrow{\sim} Fuk(\bar{X}, F) \xrightarrow{\sim} Fuk(\bar{X}, \mathfrak{c}_F).$$

\subsection{LG model}  \label{sec:lef}

Partially wrapped Fukaya categories can be used to formulate homological mirror symmetry for
Fano varieties.  For example, the mirror to $\mathbb{P}^1$ 
should be somehow associated to the function $W(z) = z + z^{-1}$ on $\mathbb{C}^*$.  
We interpret this to mean that we should form a Liouville sector from $\mathbb{C}^*$ by 
deleting the neighborhood of a fiber $W^{-1}(-\infty)$ at infinity.  
In this special case, any reasonable interpretation of the above description should result 
in the sector on the left-hand side of Figure \ref{fig:mirrorexample}.   

More generally, we would like to obtain a Liouville sector from a function $W: (\CC^*)^n \to \CC$, 
as such functions were predicted by Hori and Vafa \cite{HV} to provide mirrors to toric varieties. 
Na\"ively, one could attempt to produce a sector from this data as follows: take a half-plane $\HH \subset \CC$ containing all the critical values
(including those associated to critical points at infinity) of $W$, and take $W^{-1}(\HH)$ as the sector.  
Strictly speaking, however, $W^{-1}(\HH)$ is not generally conical at infinity for the restriction 
of the most natural Liouville structure on $(\CC^*)^n$, so some manipulation of exact
structures and use of cutoff functions 
would be necessary.  Similar issues arise in work of Seidel, see e.g. \cite[Sec. 3A]{Sei-more} and 
\cite[Sec. 19B]{Sei-book}. 
Instead, we use the tropical methods of \cite{M, A1} to show the following: 

\begin{proposition} \label{prop:canembedfiber}  
Fix a Newton polytope $\Delta^\vee \subset \ZZ^n$ and regular star subdivision $\mathcal{T}$ induced by some 
piecewise-linear function $\alpha$.  Consider the Laurent polynomial
$$W(z)  = \sum_{m\in V}t^{-\alpha(m)}z^m. $$ 
There is a real codimension-2 symplectic submanifold $F_\Sigma$ of $(\CC^*)^n$ such that: 
\begin{itemize}
\item \cite{A1} For $t \gg 0$, there is an isotopy of symplectic submanifolds between $F_\Sigma$ and a
general fiber $F_W$ of $W$.
\item There is a Liouville subdomain $D \subset (\CC^*)^n$, completing to $(\CC^*)^n$, such that 
$\partial D \cap F_\Sigma$ is a Liouville subdomain of $F_\Sigma$, completing to $F_\Sigma$. 
\end{itemize}
\end{proposition}

As indicated, the first item is proven in \cite{A1}, in a form we recall in Lemma \ref{lem:lhconstr}.  
The second item follows from our further calculations that  that the skeleton of $F_\Sigma$ is contained
in the boundary of some subdomain $D$ (Theorem \ref{thm:hypskel}), and  moreover, along the skeleton, $F_\Sigma$ 
is nowhere tangent to the Liouville vector field 
of the ambient $(\CC^*)^n$ (Lemma \ref{lem:nottangent}).  
Indeed, then we may deform slightly $\partial D$ along the Liouville field in order to contain some neighborhood of the skeleton of $F_\Sigma$.

It is the sector associated to the particular
pair $(D, F_\Sigma \cap D)$ constructed by our proof of this proposition that is used in this
article.  In particular in Theorem \ref{maintheorem}, when we assert `there is a Liouville structure on $F_W$' 
we mean that we pull back the Liouville structure mentioned above under the symplectomorphism $F_W \cong F_\Sigma$. 

Of course, we expect that any other reasonable construction of such a pair from $W$ 
will be deformation equivalent to ours, in particular giving the same Fukaya category.  

\subsection{Sheaves} 

A prototypical example of a Liouville manifold is the cotangent bundle $T^*M$ of a closed manifold without boundary; the skeleton for the
usual ``$pdq$'' form is the zero section.  If $M$ had boundary, the cotangent bundle would naturally be a Liouville
sector, again with the zero section as skeleton. 
An open set $U \subset M$ determines an inclusion
of Liouville sectors
$T^*U \subset T^*M$: the stopped boundary of $T^*U$ is the restriction of the cotangent bundle to the boundary of $U$. 
Lifting a cover of $M$ gives a cover of $T^*M$ by Liouville sectors, whose intersections are again Liouville sectors (with corners). 
The covariantly functorial \cite{GPS1} assignment $U \mapsto Fuk(T^* U)$ thus
defines a precosheaf of categories on $M$.  

Suppose we knew this precosheaf were a cosheaf.  Then we could compute its global sections from the local data. 
Indeed, the Fukaya category of the cotangent bundle of a disk 
is equivalent to the category of chain complexes, so the cosheaf in question would be a locally constant cosheaf
of categories.  Recall that the
$\infty$-groupoidal version of the Seifert-van Kampen theorem asserts that the fundamental 
higher groupoid of a space is the global sections of a locally constant cosheaf of spaces with stalk a point.  Linearizing this, we see
that a locally constant cosheaf of $A_\infty$ categories with stalk the category of chain
complexes has global sections  (a twisted version of) the category of modules 
over the algebra of chains on the based loop space of $M$.  Thus, the Fukaya category
of a cotangent bundle is the category of modules over chains on the based loop space.  This final statement 
is originally a result of Abouzaid, by a different argument \cite{A3}.  

Kontsevich's localization conjecture \cite{Kon-local} asserts that the existence of a similar cosheaf $\cF uk$
over the skeleton of any Weinstein manifold (e.g. the complement of an ample divisor in a smooth projective variety), whose global sections should recover the wrapped Fukaya category.  

A variant, which gives a local-to-global principle without any mention of skeleta, 
is the main result of \cite{GPS2}, which asserts that
the Fukaya category satisfies descent with respect to sectorial covers.\footnote{
To deduce Kontsevich's statement from \cite{GPS2}, 
one would want to know further that appropriate open covers 
of the skeleton lift to sectorial covers.  It is expected that such a lifting is not difficult to construct
in general.  In the case of relevance to this article, it is likely possible to construct such a cover by hand, 
though we will not do it here, as we do not invoke this result (instead we use \cite{GPS3}).}  This result, together with the well known calculation of Fukaya
categories of disks with stops at the boundary (for a very short calculation, see \cite[Ex. 1.22]{GPS2}), 
can easily be used to make the discussion
of Section \ref{pictures} above completely rigorous.  

In the body of this article we will need a further elaboration of Kontsevich's conjecture, 
formulated by Nadler \cite{Nwms} (and further
elaborated in \cite{S, NS}), which identifies Kontsevich's conjectural cosheaf of categories on 
the skeleton with a certain cosheaf of a combinatorial-topological nature which is constructed
directly from the microlocal sheaf theory of \cite{KS}.  This conjecture is established\footnote{Strictly 
speaking, this is established in the ``stably polarized'' case, which includes the examples of interest here.} in \cite{GPS3}, 
using the theory developed in \cite{GPS1, GPS2} and the antimicrolocalization
lemma of \cite{NS}.  

\subsection{Proof of Theorem \ref{maintheorem} } \label{sec:maindiagram}

Here we give the proof of Theorem \ref{maintheorem}, modulo the results which are 
the essential mathematical contents of the present article.  We fix the following notations
for the relevant toric data.  (A brief review of relevant algebraic geometry of toric
varieties is included in Section \ref{sec:toric}.)

Let $\TT$ be a real $n$-dimensional torus.  Let $M$ and $M^\vee$ be its lattices of characters
and cocharacters, respectively.  For an abelian group $A$, we write $M_A:=M\otimes_\ZZ A$. We can then write the torus $\TT$ and its dual as
$$
\xymatrix{
\TT=M^\vee_{\RR/\ZZ}=M^\vee_\RR/M^\vee, 
&\TT^\vee =M_{\RR/\ZZ}=M_\RR/M.
}
$$
We denote the corresponding complex tori by
$\TT_\CC=M^\vee_{\CC^\times},$ $\TT_\CC^\vee=M_{\CC^\times}.$

These complex tori are naturally tangent
bundles, $\TT_\CC^\vee\cong T\TT^\vee=(M_\RR/M)\times M_\RR$, but we will 
 choose an inner product 
to produce an identification with the cotangent bundle $T \TT^\vee \cong T^*\TT^\vee$.  We always regard the latter as an exact symplectic manifold
carrying the canonical (``$p dq$'') Liouville structure.

\begin{remark}
	The above choice of inner product is an essential feature of mirror symmetry: even in the most basic mirror pair of $\CC^\times$ and $T^*S^1,$ the inner product is necessary for constructing dual torus fibrations over a shared SYZ base. In our setting, this inner product allows us to present a hypersurface with Newton polytope $\Delta^\vee,$ naturally a complex submanifold of $\TT^\vee_\CC,$ as a symplectic submanifold of $T^*\TT^\vee.$
\end{remark}

Our mirror-symmetric setup is as follows.  
Let $\Delta^\vee\subset M_\RR^\vee$ be an integral polytope containing the origin.  
Choose a regular star-shaped triangulation of $\Delta^\vee$; equivalently, 
choose a smooth quasiprojective stacky fan $\Sigma \subset M_\RR^\vee$ whose stacky primitives lie on $\partial \Delta^\vee$ and have convex hull $\Delta^\vee$. This determines a toric stack $\bT_\Sigma$ partially compactifying $\TT_\CC$, and we denote its toric boundary by $\partial \bT_\Sigma$.  

\begin{remark}
To state results in their natural generality, we use the toric stacks of \cite{BCS}.  
For the purpose of understanding the new ideas in this paper, this can be entirely ignored. 

Very briefly, 
toric stacks are 
smooth Deligne-Mumford stacks associated to the data of a ``smooth stacky fan'' $\Sigma$, 
which is a simplicial fan together with a choice of integer point along each ray.  
We term these chosen integer points the ``stacky primitives''.  
The coarse moduli space of the toric stack is the toric variety which would ordinarily
correspond to the underlying simplicial fan.  

Even in the setting of reflexive polytopes, one must in general allow stacks to get the correct 
category of coherent sheaves for the purposes of mirror symmetry; this is due to the fact that toric varieties
do not in general admit crepant resolutions.
Of course, if we begin with a smooth fan, no discussion of toric stacks is necessary.

The added generality provided by allowing toric stacks can be seen by the following lemma: 

\begin{lemma}\label{lem:regtri}
Every convex polytope containing the origin is the convex hull of the stacky primitives of a smooth quasi-projective stacky fan.  
\end{lemma}
\begin{proof}
The quasi-projectivity condition is that the triangulation induced by the fan is {\em regular}, in the sense of being the corner locus
of a piecewise-linear function $\alpha:\Delta^\vee\to\RR$.  Choose an integer point in the polytope, and let $\alpha_0$ be the piecewise linear function which is $1$ at 
the origin, and $0$ at all facets of the boundary not containing the origin.  For each facet of the polytope, $\tau$, choose some $\alpha_\tau$ inducing a regular triangulation
of $\tau$.  Then take the function $\alpha = \alpha_0 + \sum \epsilon_\tau \alpha_\tau$ for small $\epsilon_\tau$.  
(We thank Allen Knutson for this argument.)
\end{proof}
\end{remark}

We take
$W:\TT^\vee_\CC\to\CC$ a Laurent polynomial whose Newton polytope is $\Delta^\vee$.  (How to choose
this polynomial will be discussed further below, though generic choices are isotopic and hence will
determine the same categories.)

Finally, we will need a certain conical (singular) Lagrangian $\LL_\Sigma \subset T^* \TT^\vee$ introduced in \cite{FLTZ}
to study toric mirror symmetry.    We recall its definition in Section \ref{sec:FLTZ}.  

\begin{figure} 
\begin{equation} \label{maindiagram} 
\begin{tikzcd}[column sep = tiny]
 \Coh(\partial \bT_\Sigma) \ar{r} \ar[equals]{d}{\text{Thm. \ref{thm:ccc-infty}}} & \Coh(\bT_\Sigma)  \ar[equals]{d}{\text{\cite{Ku}}} \\ 
 \mu sh_{\partial \LL_\Sigma}(\partial \LL_\Sigma)^c \ar{r} \ar[equals]{d}{\text{\cite[Thm 1.4]{GPS3}}} & 
 \Sh_{\LL_\Sigma}(\TT^\vee)^c  \ar[equals]{d}{\text{\cite[Thm 1.1]{GPS3}}}  \\ 
  Fuk(F_\Sigma) \ar{r} & Fuk(T^*\TT^\vee, F_\Sigma) 
\end{tikzcd}
\end{equation}
\caption{\label{fig:maindiagram} The commutative diagram organizing our proof of Theorem 
\ref{maintheorem}.}
\end{figure}

The proof of Theorem \ref{maintheorem} proceeds by establishing the commutative
diagram in Figure \ref{fig:maindiagram}.  Indeed, the theorem follows
from the left column (whose notation we have not yet explained in its entirety), together with the
fact that $F_W$ is deformation equivalent to a general fiber of $W$ (per 
Proposition \ref{prop:canembedfiber}) and hence has the same Fukaya category.  The 
full diagram gives a functoriality result connecting mirror symmetry for the toric 
variety and for its boundary.  In fact, we will prove even stronger functoriality results 
on our way to the theorem. 

Let us now explain the diagram in detail. We have by now introduced all the geometric players: 
the real torus $\TT$ and its dual real torus $\TT^\vee$; 
the toric variety $\bT_\Sigma$ and its boundary $\partial \bT_\Sigma$; the 
\cite{FLTZ} Lagrangian $\LL_\Sigma$ and its Legendrian boundary at infinity 
$\partial \LL_\Sigma$; the Laurent polynomial $W:\TT^\vee_\CC\to\CC$,
which, under a choice of isomorphism $T^*\TT^\vee = T \TT^\vee = \TT^\vee_\CC$, 
becomes $W:T^*\TT^\vee \to\CC$; and finally $F_\Sigma$, the deformation 
$F_W$ of a general fiber of $W$.  

For an algebraic scheme (or stack) X, we write $\Coh(X)$ for the dg category of complexes 
of sheaves with coherent cohomology on $X$, localized at quasi-isomorphisms.  The top horizontal
arrow is the pushforward.

For $\Lambda \subset T^*M$, the 
notation $Sh_\Lambda(M)$ means the category of sheaves whose microsupport
is contained in $\Lambda$.  (When $\Lambda$ is instead a Legendrian in $S^*M$, we use 
the same notation of sheaves whose microsupport at infinity is contained in $\Lambda$.) 
This notion is introduced and studied in \cite{KS}.  Following
more modern conventions, and unlike in 
\cite{KS}, by $Sh$ we mean the dg category of all complexes of sheaves
localized at the acyclic complexes, rather than the bounded derived category.  We write 
$Sh(-)^c$ for the subcategory of compact objects, {\em i.e.}, the ``wrapped microlocal sheaves'' of 
\cite{Nwms}.  

The particular example of $\Sh_{\LL_\Sigma}(\TT^\vee)$ is the subject of 
\cite{FLTZ2, Tr, Ku}.  The top right vertical equality is the main result of \cite{Ku},\footnote{
When $\Sigma$ is not smooth and proper, even the functor is new in \cite{Ku}: 
the functor described in \cite{B, FLTZ, Tr} 
takes values in quasi-coherent sheaves, and it is necessary to lift this functor to take values in ind-coherent sheaves.
} building
on \cite{FLTZ2, Tr}.   This equality holds for any $\Sigma$, without the hypotheses of
smoothness or quasiprojectivity.

For $\Lambda \subset T^*M$ or $\Lambda \subset S^*M$, the notation $\mu sh_\Lambda$ 
denotes a certain sheaf of categories on $\Lambda$ constructed out of the microlocal
sheaf theory, called the Kashiwara-Schapira stack.  We recall its properties in Section \ref{subsub:msh} below. 
For formal reasons, taking compact objects in $\mu sh_\Lambda$ gives a {\em cosheaf} of categories 
$\mu sh(-)^c$.   One of our main results is the following: 

\begin{quote} {\bf 
	For $\Sigma$ determining a smooth toric stack $\bT_\Sigma$, 
there is an isomorphism $Coh(\partial \bT_\Sigma) \cong \mu sh(\partial \LL_\Sigma)^c$ (Thm. \ref{thm:ccc-infty})
ensuring that the top square commutes.}
\end{quote}

\begin{remark}
In fact such an isomorphism exists without the smoothness hypothesis.  We do not show this here but
briefly indicate how one can, see Remarks \ref{rem: Kuwagaki} and \ref{rem: blowdown}. 
\end{remark}

As the horizontal arrows in the diagram are not fully faithful, the existence
of a morphism $Coh(\partial \bT) \to \mu sh(\partial \LL_\Sigma)$ making the top 
square commute {\em does not} imply that said morphism is an isomorphism.  
A separate argument is required. 
We then use the fact 
(explained in Section~\ref{sec: FLTZ boundary}) that $\partial \LL_\Sigma$ has a cover by mirror skeleta to the toric
varieties in $\partial \bT$, together with the fact that $Coh$ and $\mu sh$ satisfy certain local-to-global
principles, to deduce this result.   To make this work, we will need to prove a functoriality result
(``restriction is mirror to microlocalization'') for the isomorphism 
$Coh(\bT) \xrightarrow{\sim} \Sh_{\LL_\Sigma}(\TT^\vee)$.

This top square is where homological mirror symmetry happens: the sheaf categories 
are already some kind of interpretation of the $A$-model (morally: in a rescaling limit under the Liouville
flow).  

The bottom square compares the microlocal sheaf categories with the Fukaya category\footnote{This has a purpose aside 
from merely matching historical formulations: 
it is the Fukaya category
which one knows how to deform away from the large volume limit (by holomorphic disks passing through a 
compactifying boundary divisor).  However, we do not take up the study of this deformation
in the present work.} 
The engine for this is the work \cite{GPS3}, whose main results we summarize:\footnote{
To compare with what is written in \cite{GPS3}, note the canonical equivalence $Fuk(T^*M, \Lambda)^{op} = Fuk(T^*M, -\Lambda)$, 
also noted in Rem. 1.2 of that reference.}

\begin{theorem} \cite[Thm. 1.1, Thm. 1.4, Cor. 7.22]{GPS3}
Let $M$ be a real analytic manifold and $\Lambda \subset S^*M$ an isotropic subanalytic subset.  Then there 
is an equivalence of categories $Fuk(T^*M, -\Lambda) \cong Sh_\Lambda(M)^c$.  If in addition $-\Lambda$ is the 
core of a Liouville hypersurface $F$ which admits homological cocores, then there is a commutative diagram
$$
\begin{tikzcd}[column sep = tiny]
 \mu sh_{-\Lambda}(-\Lambda)^c \ar{r} \ar[equals]{d} & 
 \Sh_{-\Lambda}(M)^c  \ar[equals]{d} \\ 
  Fuk(F) \ar{r} & Fuk(T^*M, F) 
\end{tikzcd}
$$
where the top map is the left adjoint to microlocalization, the bottom map is the \cite{GPS1, GPS2} functor associated to a Liouville pair, 
and the right column is related
to the aforementioned equivalence by the canonical $Fuk(T^*M, F) \xrightarrow{\sim} Fuk(T^*M, -\Lambda)$. 
\end{theorem}

In the case at hand, our $\LL_\Sigma$ will be evidently subanalytic.  The commutative diagram asserted to exist will match
the bottom square, once we establish the following:

\begin{quote}
{\bf The construction of $F_\Sigma$ in Prop. \ref{prop:canembedfiber} may be arranged so that 
the skeleton of $F_\Sigma$ is $-\partial \LL_\Sigma$ (Thm. \ref{thm:hypskel}).}
\end{quote}

We show this by using Mikhalkin-Viro patchworking \cite{M} to deform the hypersurface 
in such a way that the calculation of the skeleton localizes to ``pairs of pants,'' where in fact
it has already been studied by Nadler \cite{Nwms}.  
Our construction will show that $\LL_\Sigma$ is the skeleton associated to a Morse-Bott
Liouville flow, hence $F_\Sigma$ admits geometric cocores (and thus homological cocores).

This completes the proof of Theorem \ref{maintheorem}, modulo 
the bolded promissory notes.   $\square$

\subsection{Other related works}

We end the introduction by attempting to situate our work in the landscape of 
homological mirror symmetry.  

Our approach has been to pass as
quickly as possible to microlocal sheaf theory, and match functorial structures
on both sides in order to reduce mirror symmetry to elementary calculations. 
Previous works in this spirit include \cite{FLTZ2, Ku, Nwms}; the particular approach
used in this article is close to what is suggested in \cite{TZ}. 
The underlying topological spaces of some of the Lagrangian skeleta 
we construct were studied earlier in \cite{RSTZ}.

Note we use the foundational work \cite{GPS1, GPS2, GPS3} 
rather than \cite{NZ, Ncs}; among other reasons, this allows us to make statements regarding
the {\em wrapped} Fukaya category.

Another strategy to approach mirror symmetry is to identify particular Lagrangians, 
compute their Floer theory, and identify the resulting algebra with some 
endomorphism algebra on the mirror.  We view this as the approach taken to 
the quartic K3 in \cite{Sei}, to toric varieties in \cite{A1, A2}, and to
and hypersurfaces in projective space in \cite{Sher1, Sher2, Sher3}. 

After finding the skeleton and corresponding cover of the hypersurface, we 
could perhaps have used \cite{A1, A2} to complete the proof of mirror symmetry for Calabi-Yau hypersurfaces.
However, this would require reworking
those arguments in the wrapped setting and establishing
the appropriate functoriality with respect to inclusion of toric divisors. 
In addition,
the works \cite{A1, A2}, as \cite{FLTZ, FLTZ2, Tr}, give only a fully faithful embedding of the coherent sheaf category
into the Fukaya category; one would need to prove generation.
In any case, the form of the results in \cite{Ku} is better adapted to our uses here.

Finally we note that in \cite{AAK}, one finds a mirror proposal for very affine hypersurfaces in
terms of a category of singularities; it is \emph{a priori} different from the category we have found here. 
The reason for the difference is that the \cite{AAK} mirrors correspond to a maximal subdivision
of $\Delta^\vee$, and we have taken a decomposition centered at a single point.  One could
try and compare algebraically the resulting categories.  For that matter, we have provided
here many mirrors, depending on the choice of point, and it should be interesting to understand
the derived equivalences between them in algebro-geometric terms.

The \cite{AAK} mirrors can also be approached directly by the methods of this paper. The main new difficulty in
carrying this out is that the amoebal complements have many bounded components, 
making it more difficult to find a contact-type hypersurface containing the skeleton.  
It is, however, possible to use a higher-dimensional version of the inductive argument in \cite{PS}. 
That proof has two essential ingredients: a gluing result and a way to move around the skeleton to allow further gluings.  
The gluing result needed is exactly our microlocalization of the theorem of Kuwagaki.  We will return elsewhere
to the question of its interaction with deformations of the skeleton.

\newpage

\section{Toric geometry} \label{sec:toric}
We recall here some standard notations and concepts from toric geometry; proofs, details, and further exposition can
be found, \emph{e.g.}, in the excellent resources \cite{F, CLS}. 

In most of this paper we will be interested in a fixed toric variety $\bT$, with dense open torus $\TT_\CC$ whose character and cocharacter lattices are denoted by $M$ and $M^\vee,$ respectively.  When we must discuss another toric variety $\bT'$, we
indicate the corresponding characters and cocharacters by $M(\bT')$ and $M^\vee(\bT'),$ respectively.  In our review here we confine ourselves
to the case of toric varieties; for toric stacks see \cite{BCS}. 

\subsection{Orbits and fans} \label{subsec:orbits}

A toric variety $\bT$ is stratified by the finitely many orbits of the torus $\TT_\CC$.  The geometry of
this stratification determines a configuration of rational polyhedral cones (the `fan') in the cocharacter space. 
We briefly review this correspondence. 

For any cocharacter $\eta: \G_m \to \TT_\CC$, one can ask whether $\lim_{t \to 0} \eta(t) \in \bT$, and if so, in which
orbit it lies.

This gives a collection of regions in $M^\vee$, and for such a region $\sigma$ we denote the corresponding orbit by 
$O(\sigma)$.  Each cone $\sigma$ is readily seen to be closed under addition; in fact, each is 
the collection of interior integral points inside a rational polyhedral cone $\sigma \subset M^\vee_\RR$.  
This collection of cones is called the {\em fan} of $\bT$.  Every face of the cone in the fan is again a cone in the fan.

A character $\chi \in M$ is by definition a map $\TT_\CC \to \G_m$, but composing with the inclusion $\G_m \to \AA^1$ determines a function
on $\TT_\CC$.  One can ask whether such a function can be extended to a given torus orbit $O(\sigma)$.   Evaluating on 
one-parameter subgroups $\eta \in \sigma$, one needs $\lim_{t \to 0} \chi(\eta(t)) = \lim_{t \to 0} t^{\langle \chi, \eta \rangle}$
to be well defined, or in other words that $\langle \chi, \eta \rangle \ge 0$.
In fact, this condition is also sufficient, and moreover the ring
of all functions on $\TT$ extending to $O(\sigma)$ is $k [ \sigma^\vee]$, where
$$\sigma^\vee = \{ \chi \in M\, | \, \langle \chi, \sigma \rangle \ge 0 \}.$$
In other words, if we write $\bT_\sigma$ for the locus in $\bT$ on which all the $k[\sigma^\vee]$ are well defined, 
the natural map $\bT_\sigma \to \mathrm{Spec}\, k[\sigma^\vee]$ is an isomorphism. 

For cones $\sigma, \tau$ in a fan, the following are equivalent: $\tau \subset \overline{\sigma}$ iff 
$\sigma^\vee \subset \tau^\vee$ iff $O(\sigma) \subset \overline{O(\tau)}$ iff $k[\sigma^\vee] \subset k[\tau^\vee]$ iff
$\bT_\tau \subset \bT_\sigma$.  As sets, 
$$\bT_\sigma = \coprod_{\tau \subset \overline{\sigma}} O(\tau)$$
$$\overline{O(\sigma)} = \coprod_{\overline{\tau} \supset \sigma} O(\tau) $$
 
\begin{defn}
Let $\Sigma$ be a fan of cones in $M_\RR^\vee.$ We denote by $\bT_\Sigma$ the toric variety determined as above by the fan $\Sigma.$
\end{defn}

\subsection{Orbit closures}
\label{sec:orbs}
Let $\sigma \subset M^\vee_\RR$ be a cone of the fan.  The corresponding orbit $O(\sigma)$ is acted on trivially by the cocharacters in $\sigma$,
hence by their span $\ZZ \sigma$.  
That is, if we write denote by $\TT_\CC/\sigma$ the complex torus $(M^\vee / \ZZ \sigma) \otimes \CC^\times$, then the $\TT_\CC$ action factors through $\TT_\CC/\sigma$.  In fact
the resulting action is free, and admits a canonical section inducing an identification $\TT_\CC/\sigma \cong O(\sigma)$.  
Note in particular that the dimension of the orbit is the codimension of the cone in the fan. 

This identification can be extended to the structure of a toric variety on the orbit closure $\overline{O(\sigma)}$.  As mentioned above, as a
set 

$$\overline{O(\sigma)} = \coprod_{\overline{\tau} \supset \sigma} O(\tau) $$

The identification of the open torus with $\TT_\CC/\sigma$ induces the following
description of the lattice of cocharacters: 

$$M^\vee(\overline{O(\sigma)}) \cong  M^\vee/\ZZ \sigma$$

The fan of $\overline{O(\sigma)}$ is obtained from the $\Sigma$ by taking the cones
$\tau$ such that $\tau \supset \overline{\sigma}$ and projecting them along 
$M^\vee \to M^\vee / \ZZ \sigma$. 

The orbit closures have the relation 
$\overline{O}_\sigma \cap \overline{O}_\tau  = \overline{O}_{\sigma \wedge \tau}$,
where $\sigma \wedge \tau$ is the smallest cone in the fan containing both $\sigma$ and $\tau$ if such a cone exists,
and by convention $\overline{O}_{\sigma \wedge \tau} = \emptyset$ if no such cone exists.  That is, 
the association $\sigma \to \overline{O}_\sigma$ is inclusion reversing.

\subsection{Fans from triangulations}\label{subsec:tri-fan}
Let $\Delta^\vee\subset M^\vee_\RR$ be an integral convex polytope containing 0.
We will be interested in stacky fans obtained from star-shaped triangulations of $\Delta^\vee.$
\begin{defn}
	A triangulation $\mathcal{T}$ of $\Delta^\vee$ is a \emph{star-shaped triangulation} if every simplex in $\mathcal{T}$ which is not contained in $\partial \Delta^\vee$ has 0 as a vertex.
\end{defn}
Such a triangulation defines a stacky fan $\Sigma$: the stacky primitives of $\Sigma$ are the 1-dimensional cones in $\mathcal{T}$, and the higher-dimensional cones in $\Sigma$ are cones on the simplices in $\mathcal{T}$ which are contained in $\partial \Delta^\vee.$

\begin{remark}
	Note that not every fan $\Sigma$ arises in the above fashion. The above construction produces only those fans $\Sigma$ satisfying the following property:
	\emph{Let $\Delta^\vee$ be the convex hull of the primitives of $\Sigma.$ Then every primitive of $\Sigma$ lies on $\partial\Delta^\vee$.}
	A more complete discussion of this restriction can be found in Section~\ref{sec:extra}.
\end{remark}

Since the subdivision $\mathcal{T}$ of $\Delta^\vee$ was a triangulation, the fan $\Sigma$ is necessarily smooth. But we would also like to require that $\Sigma$ be quasi-projective; recall that this is equivalent to the condition that the triangulation $\mathcal{T}$ be regular.
\begin{defn}
	A subdivision $\mathcal{T}$ of $\Delta^\vee$ is \emph{regular} (sometimes also called \emph{coherent}) if it is obtained by projection of finite faces of the overgraph of a convex piecewise linear function $\alpha:\Delta^\vee\cap M^\vee\to \RR.$
\end{defn}

\subsection{The toric boundary}\label{sec:toric-bdry}

In this paper, we are interested in the boundary
$\partial \bT_\Sigma$ of a toric variety $\bT_\sigma$, which by definition
is the union of the nontrivial orbit closures:
$$
\partial\bT_\Sigma=\bigcup_{0\neq \sigma\in\Sigma}\overline{O}_\sigma.
$$

In fact, we need a scheme- (or stack-)theoretic version of this statement.  Below
we always take both each $\overline{O}_\sigma$ and $\partial \bT_\Sigma$ with their 
reduced structure.  

\begin{lemma} \label{lem:boundarypushout}
In the category of algebraic stacks, $\partial \bT_\Sigma = \colim_\sigma \overline{O}_\sigma$. 
\end{lemma}
\begin{proof}
There is evidently a map $\colim_{\sigma \in {\Sigma}} \overline{O}_\sigma \to \partial \bT_\Sigma$, we must check
it is an isomorphism.  The question is \'etale local, thus we reduce to the case of 
affine toric varieties, i.e. some $\bT_\Sigma = \mathrm{Spec} \,k[\tau^\vee]$ where $\tau$ is the unique maximal cone
in the fan.  

The ring of functions $\mathcal{O}(\partial \bT_\Sigma)$ is the quotient of it by all functions which vanish
on all faces; observe that this ideal is generated by the points of the interior of $\tau^\vee$.  That is,
$\mathcal{O}(\partial \bT_\Sigma) = k[\tau^\vee]/k[\mathrm{Int}\,\tau^\vee]$.  
Meanwhile the rings of functions $\mathcal{O}(\overline{O}_\sigma)$ are the further quotients of this by 
all functions except for those on the facet of $\tau^\vee$ corresponding to $\sigma$.  

Thus we are interested in whether the map $k[\tau^\vee]/k[\mathrm{Int}\,\tau^\vee] \to \lim_{\eta < \tau^\vee} k[\eta]$ 
is an isomorphism, where the $\eta$ are the faces of $\tau^\vee$.  We can study this character by character, i.e.
separately at each integer point of $\partial \tau^\vee$.  What we must show is that the RHS is one dimensional. 
As pointed out to us by Martin Olsson, this can be seen by observing that the character $\chi$ part of the  RHS is 
computing precisely the cohomology of the normal cone to $\tau$ at the character $\chi$ -- and this cone is contractible. 
\end{proof}

We will discuss the mirror to this cover in Section~\ref{sec: FLTZ boundary}.

\section{The FLTZ skeleton}\label{sec:FLTZ}

Here we recall from \cite{FLTZ, FLTZ2, FLTZ3} the conic Lagrangian $\LL_\Sigma\subset T^*\TT^\vee$.  

\subsection{Non-stacky definition and examples} 

To a non-stacky fan $\Sigma$, 
\cite{FLTZ} associated a conic Lagrangian 
$$
\xymatrix{
\LL_\Sigma=\bigcup_{\sigma \in \Sigma}(\sigma^\perp)\times(-\sigma)\subset (M_\RR^\vee/M^\vee)\times M_\RR=T^*\TT^\vee.
}
$$
This skeleton is meant to encode the mirror geometry to the toric variety $\bT_\Sigma$, and we will term it the mirror skeleton of $\bT_\Sigma$.

We draw two examples in
Figures \ref{fig:A2} and \ref{fig:P2}.  The drawing convention is that the hairs indicate conormal directions along a hypersurface; 
likewise the circles or angles indicate conormals at a point.  Thus each picture depicts a conical Lagrangian, 
and the corresponding FLTZ skeleton is the union of this with the zero section.

\begin{example} (The mirror skeleton of $\AA^1$.)  Consider the fan in $\R$ whose sole nontrivial cone is spanned by $1 \in \R$.   
We write $\LL_1 \subset T^*S^1 = S^1 \times \RR$ for the corresponding FLTZ skeleton; it is the union of the zero-section and half a cotangent fiber at the origin: 
\[
	\LL_1 = \{(\theta,0)\mid \theta\in S^1\} \cup \{(0,\xi)\mid -\xi\in \RR_{\geq0}\}.
	\]
\end{example}

\begin{example} \label{ex: affine space skeleton} (The mirror skeleton of $\AA^n$.)  Consider the fan in $\R^n$ consisting of all cones generated by subsets of $e_1, \ldots, e_n$.  
One easily sees that the corresponding FLTZ skeleton $\LL_n \subset T^* T^n$ satisfies  $\LL_n = (\LL_1)^n$.  

Another useful description of it is as follows: 
\begin{equation}\label{eq:ind-skel}
  \LL_n = T^*_{(S^1)^n} (S^1)^n \cup \bigcup_{1\leq k \leq n} \LL_{n-1} \times (T_{S^1}^*S^1)_k,
\end{equation}
where by $T^*_{S^1}S^1$ we mean the zero section of $T^*S^1,$ with the subscript $k$ indicating that it is to be inserted in the $k$-th coordinate (with the $k,\ldots,n$ coordinates of $\LL_{n-1}$ moved forward one place).
\end{example}

\begin{figure}
\begin{center}

\begin{tikzpicture}
\fill [gray!20] (3,3) rectangle (6,6);
\draw[red, line width=1.5mm, -latex] (3, 3) -- (4, 3); 
\draw[red, thick] (3, 3) -- (6, 3); 
\draw[blue, line width=1.5mm, -latex] (3, 3) -- (3, 4); 
\draw[blue, thick] (3, 3) -- (3, 6); 
\fill[black] (0,6) circle [radius = 0.05];
\fill[black] (0,5) circle [radius = 0.05];
\fill[black] (0,4) circle [radius = 0.05];
\fill[black] (0,3) circle [radius = 0.05];
\fill[black] (0,2) circle [radius = 0.05];
\fill[black] (0,1) circle [radius = 0.05];
\fill[black] (0,0) circle [radius = 0.05];
\fill[black] (1,6) circle [radius = 0.05];
\fill[black] (1,5) circle [radius = 0.05];
\fill[black] (1,4) circle [radius = 0.05];
\fill[black] (1,3) circle [radius = 0.05];
\fill[black] (1,2) circle [radius = 0.05];
\fill[black] (1,1) circle [radius = 0.05];
\fill[black] (1,0) circle [radius = 0.05];
\fill[black] (2,6) circle [radius = 0.05];
\fill[black] (2,5) circle [radius = 0.05];
\fill[black] (2,4) circle [radius = 0.05];
\fill[black] (2,3) circle [radius = 0.05];
\fill[black] (2,2) circle [radius = 0.05];
\fill[black] (2,1) circle [radius = 0.05];
\fill[black] (2,0) circle [radius = 0.05];
\fill[black] (3,6) circle [radius = 0.05];
\fill[black] (3,5) circle [radius = 0.05];
\fill[black] (3,4) circle [radius = 0.05];
\fill[black] (3,3) circle [radius = 0.05];
\fill[black] (3,2) circle [radius = 0.05];
\fill[black] (3,1) circle [radius = 0.05];
\fill[black] (3,0) circle [radius = 0.05];
\fill[black] (4,6) circle [radius = 0.05];
\fill[black] (4,5) circle [radius = 0.05];
\fill[black] (4,4) circle [radius = 0.05];
\fill[black] (4,3) circle [radius = 0.05];
\fill[black] (4,2) circle [radius = 0.05];
\fill[black] (4,1) circle [radius = 0.05];
\fill[black] (4,0) circle [radius = 0.05];
\fill[black] (5,6) circle [radius = 0.05];
\fill[black] (5,5) circle [radius = 0.05];
\fill[black] (5,4) circle [radius = 0.05];
\fill[black] (5,3) circle [radius = 0.05];
\fill[black] (5,2) circle [radius = 0.05];
\fill[black] (5,1) circle [radius = 0.05];
\fill[black] (5,0) circle [radius = 0.05];
\fill[black] (6,6) circle [radius = 0.05];
\fill[black] (6,5) circle [radius = 0.05];
\fill[black] (6,4) circle [radius = 0.05];
\fill[black] (6,3) circle [radius = 0.05];
\fill[black] (6,2) circle [radius = 0.05];
\fill[black] (6,1) circle [radius = 0.05];
\fill[black] (6,0) circle [radius = 0.05];
\end{tikzpicture}
$\qquad \qquad \qquad$
\begin{tikzpicture}
\draw[red, thick, lefthairs] (0,0) -- (0, 6); 
\draw[red, thick, lefthairs] (6,0) -- (6, 6); 
\draw[blue, thick, righthairs] (0,0) -- (6,0); 
\draw[blue, thick, righthairs] (0,6) -- (6,6); 
\fill[gray] (0,0) -- +(-180:0.2) arc (-180:-90:0.2);
\fill[gray] (6,0) -- +(-180:0.2) arc (-180:-90:0.2);
\fill[gray] (0,6) -- +(-180:0.2) arc (-180:-90:0.2);
\fill[gray] (6,6) -- +(-180:0.2) arc (-180:-90:0.2);
\end{tikzpicture}
\end{center}
\caption{ \label{fig:A2}  The fan and FLTZ skeleton for $\AA^2$.}
\end{figure}

\begin{figure}
\begin{center}

\begin{tikzpicture}
\fill [gray!20] (0,0) rectangle (6,6);
\draw[red, line width=1.5mm, -latex] (3, 3) -- (4, 3); 
\draw[red, thick] (3, 3) -- (6, 3); 
\draw[blue, line width=1.5mm, -latex] (3, 3) -- (3, 4); 
\draw[blue, thick] (3, 3) -- (3, 6); 
\draw[green, line width=1.5mm, -latex] (3, 3) -- (2, 2); 
\draw[green, thick] (3, 3) -- (0, 0); 
\fill[black] (0,6) circle [radius = 0.05];
\fill[black] (0,5) circle [radius = 0.05];
\fill[black] (0,4) circle [radius = 0.05];
\fill[black] (0,3) circle [radius = 0.05];
\fill[black] (0,2) circle [radius = 0.05];
\fill[black] (0,1) circle [radius = 0.05];
\fill[black] (0,0) circle [radius = 0.05];
\fill[black] (1,6) circle [radius = 0.05];
\fill[black] (1,5) circle [radius = 0.05];
\fill[black] (1,4) circle [radius = 0.05];
\fill[black] (1,3) circle [radius = 0.05];
\fill[black] (1,2) circle [radius = 0.05];
\fill[black] (1,1) circle [radius = 0.05];
\fill[black] (1,0) circle [radius = 0.05];
\fill[black] (2,6) circle [radius = 0.05];
\fill[black] (2,5) circle [radius = 0.05];
\fill[black] (2,4) circle [radius = 0.05];
\fill[black] (2,3) circle [radius = 0.05];
\fill[black] (2,2) circle [radius = 0.05];
\fill[black] (2,1) circle [radius = 0.05];
\fill[black] (2,0) circle [radius = 0.05];
\fill[black] (3,6) circle [radius = 0.05];
\fill[black] (3,5) circle [radius = 0.05];
\fill[black] (3,4) circle [radius = 0.05];
\fill[black] (3,3) circle [radius = 0.05];
\fill[black] (3,2) circle [radius = 0.05];
\fill[black] (3,1) circle [radius = 0.05];
\fill[black] (3,0) circle [radius = 0.05];
\fill[black] (4,6) circle [radius = 0.05];
\fill[black] (4,5) circle [radius = 0.05];
\fill[black] (4,4) circle [radius = 0.05];
\fill[black] (4,3) circle [radius = 0.05];
\fill[black] (4,2) circle [radius = 0.05];
\fill[black] (4,1) circle [radius = 0.05];
\fill[black] (4,0) circle [radius = 0.05];
\fill[black] (5,6) circle [radius = 0.05];
\fill[black] (5,5) circle [radius = 0.05];
\fill[black] (5,4) circle [radius = 0.05];
\fill[black] (5,3) circle [radius = 0.05];
\fill[black] (5,2) circle [radius = 0.05];
\fill[black] (5,1) circle [radius = 0.05];
\fill[black] (5,0) circle [radius = 0.05];
\fill[black] (6,6) circle [radius = 0.05];
\fill[black] (6,5) circle [radius = 0.05];
\fill[black] (6,4) circle [radius = 0.05];
\fill[black] (6,3) circle [radius = 0.05];
\fill[black] (6,2) circle [radius = 0.05];
\fill[black] (6,1) circle [radius = 0.05];
\fill[black] (6,0) circle [radius = 0.05];
\end{tikzpicture}
$\qquad \qquad \qquad$
\begin{tikzpicture}
\draw[red, thick, lefthairs] (0,0) -- (0, 6); 
\draw[red, thick, lefthairs] (6,0) -- (6, 6); 
\draw[blue, thick, righthairs] (0,0) -- (6,0); 
\draw[blue, thick, righthairs] (0,6) -- (6,6); 
\draw[green, thick, righthairs] (6,0) -- (0,6);
\fill[gray!40] 
    (0,0) circle [radius=0.2];
\fill[gray!40] 
    (0,6) circle [radius=0.2];
\fill[gray!40] 
    (6,0) circle [radius=0.2];
\fill[gray!40] 
    (6,6) circle [radius=0.2];
\end{tikzpicture}
\end{center}
\caption{ \label{fig:P2}  The fan and FLTZ skeleton for $\PP^2$.} 
\end{figure}

\subsection{Stacky definition and example}

In \cite{FLTZ3}, a `stacky' version of this construction is given.
Note first that we can understand the torus $\TT^\vee$ as the Pontrjagin dual of the lattice $M^\vee:$
$$
\TT^\vee=\hatMc=\Hom(M^\vee,\RR/\BZ).
$$
Now let $\sigma\in \Sigma$ be a cone, corresponding to a face $F_\sigma$ of the polytope $\Delta^\vee.$ If $F_\sigma$ has vertices $\beta_1,\ldots,\beta_k,$ then we denote by $M^\vee_\sigma$ the quotient
$$
M^\vee_\sigma=M^\vee/\langle\beta_1,\ldots,\beta_k\rangle.
$$

Thus the group of homomorphisms $\Hom(M^\vee_\sigma,\RR/\BZ),$ which we will denote by $G_\sigma,$ is a possibly disconnected subgroup of $\hatMc=\TT^\vee.$ We  write $\Gamma_\sigma$ for the group $\pi_0(G_\sigma)$ of components of $G_\sigma.$ We  use these possibly disconnected tori to define $\LL_\Sigma$ in the general case.
\begin{defn}The \emph{FLTZ skeleton} $\LL_{{\Sigma}}\subset T^*\TT^\vee$ is the conic Lagrangian
$$
\LL_{{\Sigma}}=\bigcup_{\sigma\in\Sigma}(G_\sigma\times (-\sigma)).
$$
We will denote by $\LL_{\Sigma}^\infty$ or $\partial\LL_{\Sigma}$ the corresponding Legendrian in $T^\infty \TT^\vee$: it is the spherical projectivization of $\LL_{\Sigma}\setminus \TT^\vee.$ 
When $\Sigma$ is a non-stacky fan, this reduces to the above definition. 
\end{defn}
\begin{remark}
The relative skeleton of the Liouville sector associated to the Hori-Vafa superpotential $W$ will be $-\LL_{\Sigma}$ rather than $\LL_\Sigma$.  
This minus sign is a feature: it cancels the need for taking opposite category in the sheaf-Fukaya equivalence of \cite{GPS3}. 
\end{remark}
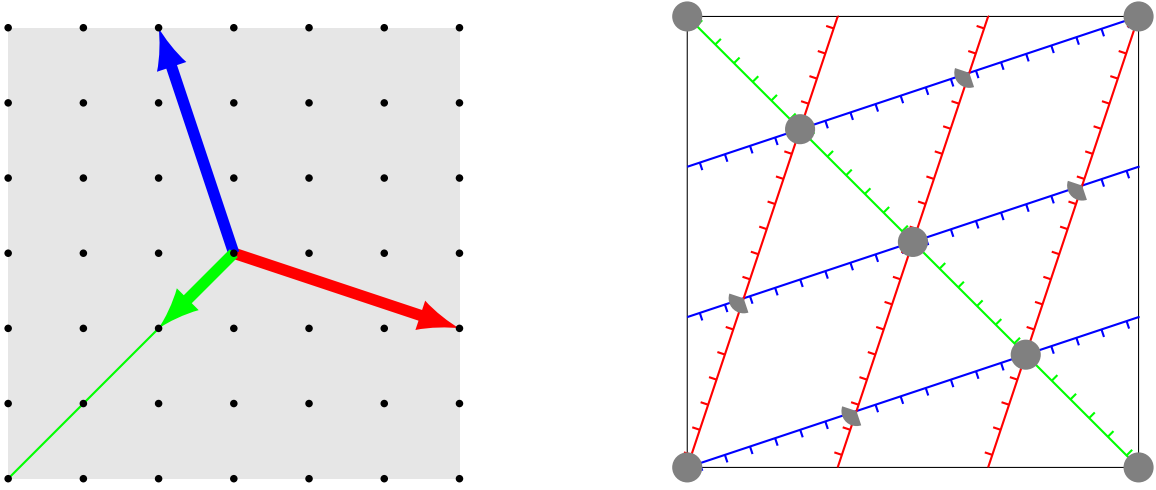
\begin{figure}[h]
\begin{center}
\begin{tikzpicture}
\fill [gray!20] (0,0) rectangle (6,6);
\draw[red, line width=1.5mm, -latex] (3, 3) -- (6, 2); 
\draw[blue, line width=1.5mm, -latex] (3, 3) -- (2, 6); 
\draw[green, line width=1.5mm, -latex] (3, 3) -- (2, 2); 
\draw[green, thick] (3, 3) -- (0, 0); 
\fill[black] (0,6) circle [radius = 0.05];
\fill[black] (0,5) circle [radius = 0.05];
\fill[black] (0,4) circle [radius = 0.05];
\fill[black] (0,3) circle [radius = 0.05];
\fill[black] (0,2) circle [radius = 0.05];
\fill[black] (0,1) circle [radius = 0.05];
\fill[black] (0,0) circle [radius = 0.05];
\fill[black] (1,6) circle [radius = 0.05];
\fill[black] (1,5) circle [radius = 0.05];
\fill[black] (1,4) circle [radius = 0.05];
\fill[black] (1,3) circle [radius = 0.05];
\fill[black] (1,2) circle [radius = 0.05];
\fill[black] (1,1) circle [radius = 0.05];
\fill[black] (1,0) circle [radius = 0.05];
\fill[black] (2,6) circle [radius = 0.05];
\fill[black] (2,5) circle [radius = 0.05];
\fill[black] (2,4) circle [radius = 0.05];
\fill[black] (2,3) circle [radius = 0.05];
\fill[black] (2,2) circle [radius = 0.05];
\fill[black] (2,1) circle [radius = 0.05];
\fill[black] (2,0) circle [radius = 0.05];
\fill[black] (3,6) circle [radius = 0.05];
\fill[black] (3,5) circle [radius = 0.05];
\fill[black] (3,4) circle [radius = 0.05];
\fill[black] (3,3) circle [radius = 0.05];
\fill[black] (3,2) circle [radius = 0.05];
\fill[black] (3,1) circle [radius = 0.05];
\fill[black] (3,0) circle [radius = 0.05];
\fill[black] (4,6) circle [radius = 0.05];
\fill[black] (4,5) circle [radius = 0.05];
\fill[black] (4,4) circle [radius = 0.05];
\fill[black] (4,3) circle [radius = 0.05];
\fill[black] (4,2) circle [radius = 0.05];
\fill[black] (4,1) circle [radius = 0.05];
\fill[black] (4,0) circle [radius = 0.05];
\fill[black] (5,6) circle [radius = 0.05];
\fill[black] (5,5) circle [radius = 0.05];
\fill[black] (5,4) circle [radius = 0.05];
\fill[black] (5,3) circle [radius = 0.05];
\fill[black] (5,2) circle [radius = 0.05];
\fill[black] (5,1) circle [radius = 0.05];
\fill[black] (5,0) circle [radius = 0.05];
\fill[black] (6,6) circle [radius = 0.05];
\fill[black] (6,5) circle [radius = 0.05];
\fill[black] (6,4) circle [radius = 0.05];
\fill[black] (6,3) circle [radius = 0.05];
\fill[black] (6,2) circle [radius = 0.05];
\fill[black] (6,1) circle [radius = 0.05];
\fill[black] (6,0) circle [radius = 0.05];
\end{tikzpicture}
$\qquad \qquad \qquad$
\begin{tikzpicture}
\draw (6,0) -- (0,0) -- (0, 6) -- (6,6) -- (6, 0); 
\draw [blue, thick, righthairs] (0,0) -- (6, 2);
\draw [blue, thick, righthairs] (0,2) -- (6, 4);
\draw [blue, thick, righthairs] (0,4) -- (6, 6);
\draw [red, thick, lefthairs] (0,0) -- (2, 6);
\draw [red, thick, lefthairs] (2,0) -- (4, 6);
\draw [red, thick, lefthairs] (4,0) -- (6, 6);
\draw [green, thick, lefthairs] (0,6) -- (6,0); 
\fill[gray] (0,0) circle [radius=0.2];
\fill[gray] (6,6) circle [radius=0.2];
\fill[gray] (0,6) circle [radius=0.2];
\fill[gray] (1.5,4.5) circle [radius=0.2];
\fill[gray] (3,3) circle [radius=0.2];
\fill[gray] (4.5,1.5) circle [radius=0.2];
\fill[gray] (6,0) circle [radius=0.2];
\fill[gray] (.75,2.25) -- +(161.565:0.2) arc (161.565:288.435:0.2);
\fill[gray] (5.25,3.75) -- +(161.565:0.2) arc (161.565:288.435:0.2);
\fill[gray] (2.25, .75) -- +(161.565:0.2) arc (161.565:288.435:0.2);
\fill[gray] (3.75, 5.25) -- +(161.565:0.2) arc (161.565:288.435:0.2);
\end{tikzpicture}
\end{center}
\caption{\label{fig:stackyskel} The stacky fan and FLTZ skeleton described in Example~\ref{ex:stackyskel}. The ``stackiness'' of this fan is due to the fact that $G_\sigma$ is nonzero for each top-dimensional cone $\sigma$; on the mirror, this is reflcted by the presence of cocircles in $\LL_\Sigma^\infty$ (the gray circles in the figure) above points other than 0.}
\end{figure}

\begin{example}\label{ex:stackyskel}
Let $\Sigma$ be the complete fan of cones in $\RR^2$ which has three
one-dimensional cones $\sigma_1,\sigma_2,\sigma_3$, spanned by the respective
vectors $(-1,3),(3,-1),$ and $(-1,-1),$ and three two-dimensional cones, which
we will denote by $\tau_{ij},$ where $\sigma_i,\sigma_j$ are the boundaries of
$\tau_{ij}.$ 

The tori $\sigma_i^\perp$ have four points of triple
intersection, and the tori $\sigma_1,\sigma_2$ have four additional points of
intersection. For any $\tau_{ij},$ the group $\Gamma_{\tau_{ij}}$ of discrete
translations of $\tau_{ij}$ is equal to the group $\sigma_i\cap\sigma_j,$ so
that for each $\tau_{ij}$ and each $p\in\sigma_i\cap\sigma_j,$ there is an
interval in the cosphere fiber $T_p^\infty\TT^\vee$ connecting the Legendrian
lifts of the tori $\sigma_i^\perp$ and $\sigma_j^\perp.$   See Figure 
\ref{fig:stackyskel}. 

The discrete data is used in the definition of the stacky skeleton to add pieces that will connect the Legendrian lifts of tori $\sigma^\perp,\tau^\perp$ in $\partial\LL_\Sigma$ over points where those tori intersect in the base $\TT^\vee.$  
\end{example}

\subsection{Recursive structure} \label{sec: FLTZ boundary} 
The Legendrian boundary of the FLTZ skeleton admits a structure that will be crucial in our proof of mirror symmetry:
it is a union of stabilized FLTZ skeleta for lower-dimensional fans, glued along their own Legendrian boundaries.
This is mirror to the fact, described above in Section~\ref{sec:toric-bdry},
that the boundary of a toric variety is the union of closures of toric orbits,
which are themselves toric varieties, as are their intersections.

Let $\Sigma$ be a (possibly stacky) fan as above. We have seen that each cone $\sigma$ in $\Sigma$
contributes a piece $G_\sigma\times (-\sigma)$ to the FLTZ Lagrangian $\LL_\Sigma.$
Write
\[
	\LL_\Sigma^\sigma:=\bigcup_{\tau\supset \sigma} G_\tau\times(-\tau)\subset \LL_\Sigma
\]
for the union, over all cones $\tau$ in which $\sigma$ is a face, of these pieces.
Observe that we have inclusion maps 
\begin{equation}\label{eq:skel-inclusion1}
	\LL_\Sigma^\tau \hookrightarrow \LL_\Sigma^\sigma
\end{equation} for any
inclusion of cones $\tau\supset\sigma.$

For a cone $\sigma,$ consider the quotient $M_\RR^\vee/\langle\sigma\rangle$
of $M_\RR^\vee$ by the subspace spanned by $\sigma$.
In this quotient, consider the reduced fan $\Sigma(\sigma)$ formed by the images of cones containing $\sigma$. 
We have seen in Section~\ref{sec:orbs} that this is the fan of the closure in
the toric variety $\bT_\Sigma$ of the toric orbit $O(\sigma)$.
We write $\LL_{\Sigma(\sigma)}$ for the FLTZ skeleton of the fan $\Sigma(\sigma),$ which we imagine as living in
the cotangent bundle of the possibly disconnected torus $T^*G_\sigma.$ (In other words, we take a disjoint
union of copies of the usual FLTZ Lagrangian for this fan in order to account for the stackiness of $\sigma.$)

Observe that for any inclusion of cones $\tau\supset \sigma,$ 
the quotient $\tau/\langle\sigma\rangle$ is the cone of conormal directions to $G_\tau$ in $G_\sigma.$ Using the factorization 
$T^*G_\sigma|_{G_\tau}=T^*G_\tau\times T^*_{G_\tau}G_\sigma,$ we can write the restriction to $G_\tau$ of the cotangent bundle to $G_\sigma$ as
\begin{equation}
  \label{eq:cotangent-prod}
  T^*G_\sigma|_{G_\tau} = T^*G_\tau \times \tau/\sigma.
      \end{equation}
      Now note that the component of $\LL_{\Sigma(\sigma)}\subset T^*G_\sigma$ contributed by $\tau$ -- the product of the perpendicular torus $G_\tau$ with the cone $\tau/\sigma$ -- is a product Lagrangian in the factorization \eqref{eq:cotangent-prod}. In other words, we have an inclusion
      \[
	G_\tau\times \tau/\sigma\hookrightarrow\LL_{\Sigma(\sigma)}.
      \]
      Moreover, any cone $\tau'$ containing $\tau$ will also contribute to $\LL_{\Sigma(\sigma)}$ a product Lagrangian contained inside \eqref{eq:cotangent-prod}; putting these all together, we get an inclusion of all of $\LL_{\Sigma(\tau)}$:
\[
	\LL_{\Sigma(\tau)}\times \tau/\sigma \hookrightarrow \LL_{\Sigma(\sigma)}.
	\]
This induces an inclusion
	\begin{equation}\label{eq:skel-inclusion2}
		\LL_{\Sigma(\tau)} \times \tau = (\LL_{\Sigma(\tau)}\times \tau/\sigma) \times \sigma\hookrightarrow
		\LL_{\Sigma(\sigma)}\times \sigma.
	\end{equation}

In particular, we may take $\sigma = 0$ and hence $\Sigma(\sigma) = \Sigma$.  Then the images of the $\LL_\tau$
agree with the aforementioned pieces:

\begin{lemma}\label{lem:FLTZ-recursive}
	The image of the map $\LL_{\Sigma(\tau)} \times \tau \hookrightarrow \LL_\Sigma$ is $\LL_\Sigma^\tau$.  
	 Moreover, under this identification, the inclusions
	\eqref{eq:skel-inclusion1} and \eqref{eq:skel-inclusion2} agree.
\end{lemma}

We can rephrase this as a statement about a cover of the Legendrian boundary $\partial \LL_\Sigma$ of the FLTZ skeleton $\LL_\Sigma.$ Let $S_\sigma\subset T^\infty\TT^\vee$ denote the boundary of $G_\sigma\times \sigma\subset T^*\TT^\vee.$
\begin{cor}\label{cor:cover1}
	The Legendrian $\partial\LL_\Sigma$ has an open cover by subsets $\Omega_\sigma\subset \partial \LL_\Sigma$, anti-indexed by the poset of nonzero cones in the fan $\Sigma,$ such that 
	$\Omega_\sigma\cong \LL_{\Sigma(\sigma)}\times S_\sigma,$ with the inclusions among these as described in Lemma~\ref{lem:FLTZ-recursive}.
\end{cor}

\subsection{T-duality description} 
In the next section we will explain how $\LL_\Sigma$ is related to the symplectic geometry
of the Hori-Vafa superpotential.  Here we informally describe another way to arrive at $\LL_\Sigma$, 
by studying the dual to the moment fibration of the toric variety.  This subsection contains no rigorous mathematical
statements and nothing in the remainder of the article depends upon it. 

Consider the example where 
$\Sigma \subset \RR$ has as cones the loci $0$,  $[0, \infty)$, and $(-\infty, 0]$, \emph{i.e.}, where $\Sigma$ is the fan whose
toric variety is the projective line $\mathbb{P}^1$.  The momentum map gives this space the structure
of a circle fibration over an interval 
whose circle fibers degenerate to zero radius at the ends.  The mirror should be again a circle fibration over an interval, 
this time with fibers degenerating to infinite radius on both ends.  Above, we made this precise 
by declaring that the mirror is the exact symplectic manifold $T^* S^1$, 
endowed with the Liouville sectorial structure in which each end of the cylinder has some stopped boundary. 
Imposing these stops results in a skeleton given by the union of the zero section and the conormal to a point.  This is precisely
the skeleton $\LL_\Sigma$ associated in \cite{FLTZ} to the fan $\Sigma.$

More generally, consider a toric Fano variety $\bT_\Sigma$, compactifying a torus $\TT$, corresponding to a fan
$\Sigma$ in $M^\vee_\RR$.  Let $\bT_\Sigma \to \Delta \subset M_\RR$ be the 
anticanonical momentum map.  The polytope $\Delta$ has the property that the cone over its polar dual $\Delta^\vee$
is just $\Sigma$. 

To find the mirror, we should take the dual torus $\TT^\vee$ as a fiber of the dual 
fibration over the polytope $\Delta^\vee \subset M^\vee_\RR$.  
This polytope will not be used to define another toric variety but rather,
under the principle that the T-dual of a collapsing fibration is a blowing up one, 
we use this polytope to define stopping conditions. 
Before, the torus spanned by the cocharacters of $\sigma$ would degenerate to radius zero
along the corresponding face; now, we want it to be impossible to go all the way around
the dualized version of this torus.  Correspondingly, for each cone $\sigma \in \Sigma$, we introduce the stop $\sigma^\perp$ over
the face of $\Delta^\vee$ whose cone is $\sigma$.  The result (up to a sign) is the skeleton $\LL_\Sigma$. 

Another derivation of $\LL_\Sigma$ by this sort of T-duality reasoning can found in \cite{FLTZ}.

\section{Pants}

\subsection{Pants} 
\label{subsec:pants}
By an $(n-1)$-dimensional pants, we mean the complement in $(\CC^\times)^{n-1}$ of a 
linear hypersurface transverse to all coordinate subspaces, or equivalently such a linear
hypersurface inside $(\CC^\times)^n$. 

Throughout our discussion of hypersurfaces in $(\CC^\times)^n,$ we wi use the map
\begin{eqnarray*}
\Log:(\CC^\times)^n & \to & \RR^n, \\
(z_1, \ldots, z_n) & \mapsto & (\log |z_1|, \ldots, \log |z_n|),
\end{eqnarray*} 
the moment map for the self-action of $(\CC^\times)^n.$

\begin{defn}
For $n\geq 1,$ the \emph{standard $(n-1)$-dimensional pants} is
$$
\xymatrix{\calP_{n-1}}=\{z_1+\cdots+z_n-1=0\}\subset(\CC^\times)^n.
$$
	The {\em amoeba} of $\calP_{n-1}$ is its image $\calA_{n-1}:=\Log(\calP_{n-1})$ in $\RR^n$ under the $\Log$ map.
\end{defn}
%
\begin{remark}
	The pants $\calP_{n-1}$ has an obvious action of the symmetric group $\Sigma_n,$ but in fact this action extends to an action of the symmetric group $\Sigma_{n+1}.$ This can be seen by writing $(\CC^\times)^n$ as the dense torus in $\PP^n,$ hence embedding $\calP_{n-1}$ as an open subset of the hypersurface $\overline{\calP}_{n-1}$ in $\PP^n$ defined by the equation
	\[\overline{\calP}_{n-1} = \{z_1+\cdots + z_n +z_{n+1}=0\}\subset \PP^n.\]
	This closed hypersurface has a manifest $\Sigma_{n+1}$ action which respects the open part $\calP_{n-1}.$ In our original coordinates, this action is  generated from the $\Sigma_n$ action by the extra generator
	\begin{equation}\label{eq:symmetry}
	(z_1,\ldots,z_n) = [z_1:\cdots:-z_{n+1}] \mapsto [-z_{n+1}:z_1:\cdots : z_n] = \left(\frac{-1}{z_n}, \frac{z_1}{z_n},\ldots,\frac{z_{n-1}}{z_n}\right),
      \end{equation}
      and the $\Log$ map becomes equivariant for the $\Sigma_{n+1}$ action on $\RR^n$ obtained by descending the symmetry \eqref{eq:symmetry} to $\RR^n$ in the evident way:
      \[(x_1,\ldots,x_n)\mapsto (-x_n,x_1-x_n,\ldots,x_{n-1}-x_n).\]

\end{remark}

Let $\Delta^\vee_{n-1}\subset\RR^{n}$ be the standard $n$-simplex, {\em i.e.}, the convex hull of the origin and standard basis vectors
$\{e_1, \ldots, e_n\}$.  
Let $\Pi_{n-1}$ be the union of positive-codimensional cones in the fan generated by $\{-e_1, \ldots, -e_n, \sum e_i \}$.  
Then $\Pi_{n-1}$ is a translate of the dual complex of $\Delta^\vee_{n-1}$, and a deformation retract of the amoeba $\calA_{n-1}$.  
The relationships between $P_{n-1}, \Delta^\vee_{n-1}, \calA_{n-1}, \Pi_{n-1}$ are the simplest instances of the general
relationship between very affine hypersurfaces and their tropicalizations, as will be recalled in detail in Section~\ref{sec:tropical-rec}.

More generally we will consider, for  $\ell_1,\ldots,\ell_n\gg0$, the {\em translated pants} 
\begin{equation}\label{eq:trans-pants}
  \calP_{n-1}^\ell = \{e^{-\ell_1}z_1+\cdots+e^{-\ell_n}z_n - 1=0\}\subset (\CC^\times)^n,
\end{equation}
	whose amoeba we denote by $\calA_{n-1}^\ell:= \Log(\calP_{n-1}^\ell).$ This amoeba can be obtained as a translation of $\calA_{n-1}$ by the vector $\ell\in\RR^n,$ which pushes it far into the first orthant.

Because the coefficients are all real, we have: 

\begin{lemma}
The components of $\partial\calA_{n-1}^\ell$ are the images of certain components of the real points of $\calP_{n-1}^\ell$. 
In particular, the component of  $\partial\calA_{n-1}^\ell$ bounding the region containing all sufficiently negative points
(which corresponds to 
the vertex 0 of the simplex $\Delta^\vee_{n-1}$) 
is the image of the real positive points of $\calP_{n-1}^\ell.$ 
\end{lemma} 
	\begin{proof}
		That the critical points of $\Log|_{\calP_{n-1}^\ell}$ are precisely the real points of $\calP_{n-1}^\ell$ is proved in \cite[Proposition 4.4]{M}. The critical values of this map certainly include the boundary components of the amoeba, and one can check that the ``bottom-left'' boundary component contains the image of the real positive points by observing that it contains the real positive point $(\frac{e^{\ell_1}}{n},\ldots,\frac{e^{\ell_n}}{n}).$
	\end{proof}

We will also want to consider certain other hypersurfaces which are naturally unramified covers of pants, or products of these with copies of $\CC^\times$.  

\begin{defn} \label{def: octopants} 
Given a map on character lattices
$T^\vee: \ZZ^k \to \ZZ^n$, consider the dual map of tori $f_T: (\CC^\times)^n \to (\CC^\times)^k$.  We write 
$\calP_T := f_T^{-1}(\calP_{k-1})$ for the variety obtained from the pants $\calP_{k-1}$ by pullback along the map $f_A$, 
and $\calA_T := \Log(\calP_T)$ for its amoeba. 

	As this variety depends only on the simplex $P=T^\vee(\Delta_{k-1}),$ we will also denote it by $\calP_P$ and its amoeba by $\calA_P$
	(where typically $T^\vee(\Delta_{k-1})$ has been named while $T$ has not); in this case, we will refer to it as the {\em $P$-pants}.
	As in equation \eqref{eq:trans-pants} above, we may also scale the coefficients of $\calP_P$ by $e^{-\ell_i}$ in order to obtain 
	the {\em translated $P$-pants} $\calP_{P}^\ell,$ whose amoeba $\calA_P^\ell$ is related to the amoeba $\calA_P$ by translation into the first orthant.
\end{defn}

We have the following relationship between amoeba:

\begin{lemma} \label{lem: amoeba cover} 
Let $T: \RR^n \to \RR^k$ be the dual of $T^\vee \otimes \RR$.  Then 
$T(\calA_T) = \calA_{k-1}$. 
\end{lemma} 

	Note that if $k < n$ and $T$ is unimodular ({\em e.g.}, an inclusion of a coordinate subspace), then 
$\calP_T \cong P_{k-1} \times (\CC^\times)^{n-k}$.

\subsection{Tailoring}\label{sec:tailoring}

\begin{prop}[{\cite[Section 6.6]{M}, \cite[Propositions 4.2, 4.9]{A1}}]\label{prop:locfact}
	Fix $\epsilon,K\in \RR$ with $0<\epsilon\ll K.$ There is a $\Sigma_{n+1}$-equivariant symplectic isotopy from $\calP_{n-1}$ to a hypersurface $\tilcalP_{n-1}$ with the following properties:
	\begin{enumerate}
		\item On the region 
		\[L_1 = \{(z_1,\ldots,z_n)\in\tcP_{n-1}\mid \log|z_1|<-K\},\]
			there is an equality \[L_1 = \{z_1\in \CC^\times\mid \log|z_1|<-K\}\times\tcP_{n-2},\]
			and analogous equalities hold on the other $n$ ends of $\tcP_{n-1}.$ 
	\item Let $L_1^\epsilon = \{(z_1,\ldots,z_n)\in \calP_{n-1}\mid \log|z_i|<-K+\epsilon\},$ and similarly for the other $n$ ends of $\calP_{n-1}.$ Then the isotopy is constant outside of $\bigcup_{i=1}^{n+1}L_i^\epsilon$.
	\end{enumerate}
	In particular, the amoeba $\tcA_{n-1} := \Log(\tilcalP_{n-1})$ differs from $\Pi_{n-1}$ only in a neighborhood of the singularities of the latter.  
	(See Figure \ref{fig: tailored} for the case $n=2$.)
\end{prop}

In Remark~\ref{rem:isodef} below, we recall from \cite[Section 4]{A1} the construction of  this isotopy, in the context of an arbitrary Newton polytope.  

\begin{defn}\label{def:legs}
	We call the regions $L_i$ defined above the {\em legs} of the pants $\tcP_{n-1}.$
\end{defn}

\begin{defn}\label{defs:tailored}
  Following \cite{Nwms}, we will call the hypersurface $\tilcalP_{n-1}$ the {\em tailored pants.} (In \cite{M}, it was called the ``localized pants.'') 
We analogously write  $\tcP_{n-1}^\ell$ and $\tcA_{n-1}^{\ell}$ for the corresponding construction applied to the translated pants $\calP_{n-1}^\ell$. 

Likewise, in the situation of Definition \ref{def: octopants}, we have a {\em tailored $P$-pants} $\tilcalP_P := f_T^{-1}(\tcP_{k-1})$ defined as the preimage of the tailored pants under the map 
$f_T$ corresponding to a choice of $k$-simplex $P=T^\vee(\Delta_{k-1}),$ and the {\em translated tailored $P$-pants} $\tcP_P^\ell$ obtained by rescaling its coefficients.
\end{defn}
Thanks to its presentation as an unramified cover of the standard tailored pants $\tcP_{n-1}$, the tailored $P$-pants $\tcP_P$ is easy to understand in terms of the tailoring construction we have already discussed. 
In particular, the analogue of Lemma \ref{lem: amoeba cover} holds for the tailored $P$-pants: 
\begin{lemma}
  $T(\Log(\tilcalP_P)) = \tcA_{n-1}$. 
\end{lemma}
The $P$-pants $\tcP_P$ also inherits from $\tcP_{n-1}$ an inductive structure on its legs, which we summarize as follows:

\begin{defn}
	The \emph{$i$th leg} $L_i$ of the $P$-pants $\tcP_P$ is the preimage, under the map $\bar{f}_A$, of the $i$th leg of the standard pants $\tcP_{n-1}.$ It is isomorphic to $\tcP_{F_i}$, where $F_i=\Conv(0,v_1,\ldots,v_{\hat{i}},\ldots,v_k)$ is the corresponding facet of $P$.
\end{defn}

\begin{figure}
\includegraphics[width=10cm]{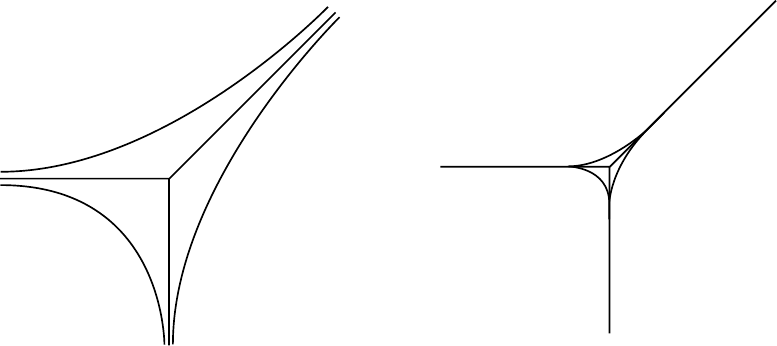}
	\caption{\label{fig: tailored}The spine $\Pi_{1},$ included in the amoebae of $\calP_{1}$ and $\tilcalP_{1}.$ (In the second picture, we have also rescaled the base of the Log map to make the situation clearer.)}
\end{figure}

\subsection{Skeleta of pants}\label{sec:skeleta2}

\subsubsection{The skeleton of $\tcP_{n-1}^\ell$}

We equip $\tcP_{n-1}^\ell$ with the restriction $\lambda$ of the symplectic primitive from $(\CC^*)^n$.  This is compatible with the recursive structure from
Proposition \ref{prop:locfact} (1): 

\begin{lemma} \label{lemma: liouville legs} 
Consider the leg $L_i$ of $\tcP_{n-1}^\ell$ for $1\leq i\leq n$.  There is an isomorphism of Liouville manifolds $L_i  \cong \tcP_{n-2}^{\ell_{\hat{i}}}\times Cyl_i$, where 
$Cyl \subset \CC^\times$ is a half-cylinder disjoint from the zero section.  The subscript $i$ on the second factor indicates that it is placed as the $i$th coordinate, 
and we write $\ell_{\hat{i}}$ for $(\ell_1,\ldots,\hat{\ell}_i,\ldots,\ell_n).$ 
\end{lemma} 
\begin{cor}
	The Liouville flow for $\lambda$ on $\tcP_{n-1}^\ell$ is complete; {\em i.e.}, $(\tcP_{n-1}^\ell, \lambda)$ is a Liouville manifold.  
\end{cor}
\begin{proof}
Recall that the product of Liouville manifolds is Liouville.  Now Lemma \ref{lemma: liouville legs} inductively characterizes the Liouville flow in the complement of a compact set.
\end{proof}

\begin{remark}
Because the original $\calP_{n-1}$ was algebraic and hence in particular a Stein submanifold of $(\CC^\times)^n$, and because the Liouville form on  $(\CC^\times)^n$
arises from a K\"ahler potential (namely $|\Log|^2$), it is also the case that the restriction of the ambient Liouville form to $\calP_{n-1}$ gives a Liouville structure on $\calP_{n-1}$.  It is presumably
true that the tailoring isotopy (recalled in Remark~\ref{rem:isodef} below from \cite[Section 4]{A1}) is an isotopy of Liouville manifolds, but we do not prove this
	here.  
\end{remark}

Recall that we write $\LL_n$ for the FLTZ skeleton mirror to affine $n$-space, as described in Example \ref{ex: affine space skeleton}.

\begin{figure}
\includegraphics[width=5cm]{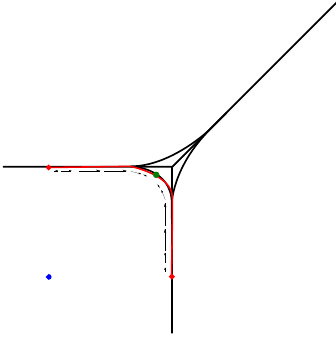}
	\caption{\label{fig:lag-simplex}The simplex $\oS_+,$ drawn in red on the amoeba of $\tcP_1,$ with its barycenter illustrated in green. Note that the vertices of $\oS_+$ are the closest points on their respective legs to the origin (blue). The arrows indicate Liouville flow along $\oS_+.$}
\end{figure}

\begin{theorem}[\cite{Nwms}]\label{thm:std-pants}
Let $\partial^0\tcA_{n-1}^\ell  \subset \tcA_{n-1}^\ell$ be the component which bounds the region of $\RR^n$ containing the all-negative orthant. 
Let $C = \Log^{-1}(\partial^0\tcA_{n-1}^\ell)\subset (\CC^\times)^n$. 

Then $C$ is a contact hypersurface, and the skeleton of $\tcP_{n-1}^\ell$ is $C\cap(-\LL_n)$ 
\end{theorem}

\begin{proof}
We proceed by induction on the dimension of the pants, the case $n=1$ being trivial.  
Many of the ideas of the proof can be seen in 
the illustration of Figure \ref{fig:lag-simplex}.

Let us consider the legs of $\tcP_{n-1}^\ell$. 
From Lemma \ref{lemma: liouville legs}, it is clear that any zero of the Liouville vector field contained in the leg $L_i$ 
must be contained inside the zero-section, {\em i.e.}, the unit circle, of its $\CC^\times_i$ factor; in other words, any
	zero of the Liouville vector field on $L_i$ must project under the $\Log$ map to the $i$th coordinate hyperplane in $\RR^n.$
	In particular, no vanishing happens on the $(n+1)$th leg of the pants, since any vanishing must be contained in the hyperplane
	given by the sum of the coordinate directions, and the translation by $\ell$ ensures that this hyperplane is
	disjoint from the leg $L_{n+1}$. 

Moreover, the preimage in $\tcP_{n-1}^\ell$ of the coordinate hyperplanes in $\RR^n$ is entirely contained in the legs, 
and stable under the Liouville flow.    By Lemma \ref{lemma: liouville legs} and the induction
hypothesis, the portion of the skeleton contained in $L_i$ is $(-\LL_{n-1}\times (T^*_{S^1}S^1)_i)\cap C$, 
using the notation of Equation (\ref{eq:ind-skel}) of Example \ref{ex: affine space skeleton}.   
By comparing that equation to the statement of this theorem, we see that our remaining task is to show 
there is exactly one more component of the skeleton, and to identify it with the intersection
	of $C$ with the positive real points of $\tcP_{n-1}^\ell$.

	Away from the legs of the pants $\tcP_{n-1}^\ell$, the map $\Log$ is a local diffeomorphism everywhere except the real points $R:= \tcP_{n-1}^\ell\cap \RR^n$. 
Let $z=(z_1,\ldots,z_n)\in R$ be a real point where the Liouville vector field vanishes.   
The equation of the pants $\sum e^{- \ell_i} z_i = 1$ prevents all $z_i$ from being negative; if 
$z_i$ is positive and $z_j$ is negative, 
then the Liouville vector field will point partially in the direction of $(z_1,\ldots,z_i+\epsilon,\ldots,z_j-\epsilon,\ldots,z_n)$ and in particular will be nonzero at $z$. 
Thus $z \in R_+ = \tcP_{n-1}^\ell\cap (\RR_{>0})^n$. In order that $z$ not lie in the legs, it must be contained in
	\[
		S_+ = \{z\in R_+\mid \Log(z)\in(\RR_{> 0})^n\}.
		\]

Recall that $\Log$ restricts to a diffeomorphism $R_+\to \partial^0\tcA^\ell_{n-1}$ from $R_+$ to the inner boundary component of the tailored amoeba.

Since $S_+$ is contained inside the real points of $\tcP_{n-1}^\ell$, the Liouville form vanishes on its tangent vectors, so it is preserved by the Liouville vector field. 
The Liouville flow increases distance to $0\in \RR^n$ under the Log projection, and the embedding of $\Log(S_+)$ in $\RR^n$ is concave and symmetric under exchange of coordinates.  
Hence the Liouville field everywhere points along $S_+$ toward the barycenter of $S_+$.  This barycenter gives the sole
remaining zero of the Liouville form, and it contributes its stable cell $S_+$ to the skeleton.  
\end{proof}

\begin{remark}
The closure  $\overline{S}_+$ of the region $S_+$ is an $(n-1)$ simplex, each facet of which is contained in one of the legs, and whose boundary projects 
	to the intersection of the amoeba with the coordinate hyperplanes. The case $n=2$ is depicted in Figure \ref{fig:lag-simplex}.
\end{remark}

\subsubsection{Skeleta for $P$-pants}
Let $P=\Conv(0,v_1,\ldots,v_k)\subset M_\RR^\vee$ be a simplex. In Definition~\ref{def: octopants}, we described the {\em $P$-pants} $\tilcalP_P\subset\TT^\vee_\CC$ obtained as a cover of the pants $\calP_{n-1}$, and in Definition \ref{defs:tailored} we described its tailored translated version $\tcP_{n-1}$. 

After choosing an inner product on $M_\RR$ and hence respective symplectic and Liouville forms $\omega$ and $\lambda$ on $\TT^\vee_\CC\cong T^*\TT^\vee$, we can restrict these to the translated tailored $P$-pants $\tcP_P^\ell$ to equip this space with the structure of a Liouville manifold. As for the standard pants, we will be interested in computing the Lagrangian skeleton of $(\tcP_P^\ell,\lambda),$ closely following the calculation in Theorem~\ref{thm:std-pants}.

let $\Sigma_P$ be the stacky fan whose primitives are the nonzero vertices of $P$.
As in the statement of Theorem~\ref{thm:std-pants}, let $\partial^0\calA_P^\ell$ be the component of the amoeba boundary $\partial \calA_P^\ell$ bounding the ``lower-left'' orthant of $\RR^n$. We will be interested in the contact hypersurface $C_P\subset \TT_\CC^\vee$ lying above this boundary:
$$
C_P:=\{z\in \TT_\CC^\vee\mid \Log(z)\in \calA_P^0\}.
$$

As in Section~\ref{sec:FLTZ}, let $G_\sigma$ be the possibly disconnected torus $\Hom(M^\vee_\sigma,\RR/\BZ),$ where $M^\vee_\sigma$ is the quotient of $M^\vee$ by the vertices of the stacky primitives in $\sigma.$ This defines a Lagrangian
\[
	\LL_{\Sigma_P}:=\bigcup_{\sigma\in\Sigma_P} G_\sigma\times\sigma\subset T^*\TT^\vee;
	\]
using the inner product, we can treat $\Sigma$ as a fan of cones in $M_\RR$ and hence $\LL_{\Sigma_P}$ as a subset of $\TT^\vee_\CC.$

For $1\leq i\leq k,$ write $\Sigma_P^i$ for the fan of cones on the $(k-1)$ vectors $v_1,\ldots,v_{\hat{i}},\ldots,v_k.$ As was the case for the standard pants, we find it helpful to rewrite the FLTZ Lagrangian as a union
\[
	\LL_{\Sigma_P} = (G_{\Sigma_P}\times\Sigma_P)\bigcup_{1\leq i \leq k}\LL_{\Sigma_P^i}
	\]
	of one new piece (where we write $\Sigma_P$ for the big cone in the fan), living in the cotangent fibers over the points $G_{\Sigma_P},$ and FLTZ skeleta for lower-dimensional cones of $\Sigma$.

\begin{lemma}\label{lem:nonstd}
There is an equality
$$
	\Lambda_P=C_P\cap (-\LL_{\Sigma_P})
$$
between the skeleton $\Lambda_P$ of $\tilcalP^\ell_P$ and 
the intersection of the contact hypersurface $C_P$ with the negative stacky FLTZ Lagrangian for $\Sigma_P.$
\end{lemma}
\begin{proof}
	The proof of Theorem~\ref{thm:std-pants} proceeded by induction on dimension, using the fact that each leg $L_i$ of the stnadard pants was itself (the product of $\CC^\times$ with) a pants one dimension lower. The proof here follows the same strategy: 
	we need to consider here $P$-pants for all $P$ (not necessarily top-dimensional), but as before we induct on the dimension of $P$.

	For clarity, we spell out explicitly the base case, when $P=\Conv(0,v)$ is 1-dimensional. In this case, 
	the tailoring construction is unnecessary, since
	$\calP_P^\ell$ is the hypersurface defined (in coordinates $z=(z_1,\ldots,z_n)$) by $\{z_1^{v_1}\ldots z_n^{v_n} = e^\ell\},$ whose amoeba $\Log(\calP_P^\ell)$ is the hyperplane
	\[
		\calA_P^\ell = \{v_1x_1+\cdots v_nx_n = \ell\}\subset \RR^n\cong M_\RR.
	\]
	In other words, the hypersurface $\calP_P^\ell$ is a copy of $(\CC^\times)^{n-1}$, with its symplectic and Liouville form restricted from those of the ambient $(\CC^\times)^n.$ Hence its Liouville vector field is given by the gradient of the restriction of the Morse-Bott function $|\Log|^2.$
	The critical locus of this function is the
	fiber of $\calP_P^\ell$ over the point $p\in \calA_P^\ell$ nearest to $0\in M_\RR,$ which is a manifold of minima for $|\Log|^2|_{\calP_P^\ell}.$ 
	As $v$ is the normal vector to the hyperplane $\calA_P^\ell,$
	the point $p$ is the point of $\tcA_P^\ell$ where it intersects the ray defined by $v.$ The fiber over this point is the preimage, under the covering map $f_A,$ of the corresponding fiber of the standard pants: this is the subtorus $G_v\subset \TT^\vee.$

	We now assume by induction that we have proven the lemma for all $P'$-pants with $\dim(P')<n,$ and we return to the case where $P=\Conv (0,v_1,\ldots v_n)$ is an $n$-simplex. From this point the proof follows very closely the proof for the standard pants. As in that case, we first investigate the legs $L_1,\ldots,L_n$ of $\tcP_P^\ell.$ Each of these is itself a $P'$-pants, for $P' = \Conv(0,v_1,\ldots,v_{\hat{i}},\ldots,v_n),$ and by induction we know that the vanishing of the Liouville vector field on leg $L_i$ contributes to the skeleton of $\tcP_P^\ell$ the piece $\LL_{\Sigma_P^i}\cap C_P$. It remains for us to determine the vanishing loci of the Liouville vector field on the interior of the pants. (As for the standard pants, it is obvious that no vanishing happens on the final leg.)

	We now consider the simplex $\oS_+ =\Ain_P\cap \Sigma_P,$ where we write $\Sigma_P$ for the top-dimensional cone in the fan, and we write $p\in \oS_+$ for the point in the interior of $\oS_+$ which is closest to 0. Let $S_+$ denote the preimage of the interior of this simplex, which is now a disjoint union $S_+ = \bigsqcup_{i=1}^d S_+^i$ of $d=\on{vol}(P)$ open simplices $S_+^i.$ Each of these simplices is preserved by the Liouville flow, which flows each simplex to the point lying over $p$, on which the Liouville field vanishes. Hence the remaining pieces of the skeleton are the open simplices $S^i_+,$ each of which is mapped diffeomorphically by $\Log$ onto the interior of $\oS_+.$ As $\oS_+$ is the intersection of $\Ain_P$ with the big cone in $\Sigma_P,$ and the fiber in $\tcP_P^\ell$ over a point in $\oS_+$ is the discrete group $G_\Sigma$, this is the desired extra piece in $\LL_{\Sigma_P}\cap C_P.$

	Finally, if there were any other vanishing of the Liouville form in the interior of the pants, it would have to lie over a critical value of $\Log$. These critical values are just the preimage (under the cover $f_A$) of the real points of $\tcP_{n-1},$ and we have already seen in the proof of Theorem~\ref{thm:std-pants} that the Liouville vector field is nonvanishing there.
\end{proof}
A crucial point is that the above result holds in the case of a simplex $P$ with arbitrary volume,
obtained as a cover of the standard simplex $\Delta_{n-1}.$  For instance, when $n=2$, the $P$-pants $\tilcalP_P$ may be higher genus. 

\begin{example}
Let $\Delta\subset\RR^2$ be the simplex with vertices $\{(0,0),(2,0),(0,2)\},$ so that the corresponding stacky fan $\Sigma$ is a stacky fan for the stack $\AA^2/(\ZZ/2\times\ZZ/2).$  We draw the stacky fan and FLTZ skeleton in Figure \ref{fig:stackyfanA2mod22}. The boundary $\partial\AA^2/(\ZZ/2\times\ZZ/2)$ matches the mirror skeleton pictured in Figure~\ref{fig:stacky-cover}.

\begin{figure}[h]
\begin{center}
\begin{tikzpicture}
\fill [gray!20] (3,3) rectangle (6,6);
\draw[red, line width=1.5mm, -latex] (3, 3) -- (5, 3); 
\draw[red, thick] (3, 3) -- (6, 3); 
\draw[blue, line width=1.5mm, -latex] (3, 3) -- (3, 5); 
\draw[blue, thick] (3, 3) -- (3, 6); 
\fill[black] (0,6) circle [radius = 0.05];
\fill[black] (0,5) circle [radius = 0.05];
\fill[black] (0,4) circle [radius = 0.05];
\fill[black] (0,3) circle [radius = 0.05];
\fill[black] (0,2) circle [radius = 0.05];
\fill[black] (0,1) circle [radius = 0.05];
\fill[black] (0,0) circle [radius = 0.05];
\fill[black] (1,6) circle [radius = 0.05];
\fill[black] (1,5) circle [radius = 0.05];
\fill[black] (1,4) circle [radius = 0.05];
\fill[black] (1,3) circle [radius = 0.05];
\fill[black] (1,2) circle [radius = 0.05];
\fill[black] (1,1) circle [radius = 0.05];
\fill[black] (1,0) circle [radius = 0.05];
\fill[black] (2,6) circle [radius = 0.05];
\fill[black] (2,5) circle [radius = 0.05];
\fill[black] (2,4) circle [radius = 0.05];
\fill[black] (2,3) circle [radius = 0.05];
\fill[black] (2,2) circle [radius = 0.05];
\fill[black] (2,1) circle [radius = 0.05];
\fill[black] (2,0) circle [radius = 0.05];
\fill[black] (3,6) circle [radius = 0.05];
\fill[black] (3,5) circle [radius = 0.05];
\fill[black] (3,4) circle [radius = 0.05];
\fill[black] (3,3) circle [radius = 0.05];
\fill[black] (3,2) circle [radius = 0.05];
\fill[black] (3,1) circle [radius = 0.05];
\fill[black] (3,0) circle [radius = 0.05];
\fill[black] (4,6) circle [radius = 0.05];
\fill[black] (4,5) circle [radius = 0.05];
\fill[black] (4,4) circle [radius = 0.05];
\fill[black] (4,3) circle [radius = 0.05];
\fill[black] (4,2) circle [radius = 0.05];
\fill[black] (4,1) circle [radius = 0.05];
\fill[black] (4,0) circle [radius = 0.05];
\fill[black] (5,6) circle [radius = 0.05];
\fill[black] (5,5) circle [radius = 0.05];
\fill[black] (5,4) circle [radius = 0.05];
\fill[black] (5,3) circle [radius = 0.05];
\fill[black] (5,2) circle [radius = 0.05];
\fill[black] (5,1) circle [radius = 0.05];
\fill[black] (5,0) circle [radius = 0.05];
\fill[black] (6,6) circle [radius = 0.05];
\fill[black] (6,5) circle [radius = 0.05];
\fill[black] (6,4) circle [radius = 0.05];
\fill[black] (6,3) circle [radius = 0.05];
\fill[black] (6,2) circle [radius = 0.05];
\fill[black] (6,1) circle [radius = 0.05];
\fill[black] (6,0) circle [radius = 0.05];
\end{tikzpicture}
$\qquad \qquad \qquad$
\begin{tikzpicture}
\draw[red, thick, lefthairs] (0,0) -- (0, 6); 
\draw[red, thick, lefthairs] (3,0) -- (3, 6); 
\draw[red, thick, lefthairs] (6,0) -- (6, 6); 
\draw[blue, thick, righthairs] (0,0) -- (6,0); 
\draw[blue, thick, righthairs] (0,3) -- (6,3); 
\draw[blue, thick, righthairs] (0,6) -- (6,6); 
\fill[gray] (0,0) -- +(-180:0.2) arc (-180:-90:0.2);
\fill[gray] (3,0) -- +(-180:0.2) arc (-180:-90:0.2);
\fill[gray] (6,0) -- +(-180:0.2) arc (-180:-90:0.2);
\fill[gray] (0,3) -- +(-180:0.2) arc (-180:-90:0.2);
\fill[gray] (3,3) -- +(-180:0.2) arc (-180:-90:0.2);
\fill[gray] (6,3) -- +(-180:0.2) arc (-180:-90:0.2);
\fill[gray] (0,6) -- +(-180:0.2) arc (-180:-90:0.2);
\fill[gray] (3,6) -- +(-180:0.2) arc (-180:-90:0.2);
\fill[gray] (6,6) -- +(-180:0.2) arc (-180:-90:0.2);
\end{tikzpicture}
\caption{\label{fig:stackyfanA2mod22} The stacky fan and FLTZ skeleton for $\AA^2/(\ZZ/2 \times \ZZ/2)$.}
\end{center}
\end{figure}
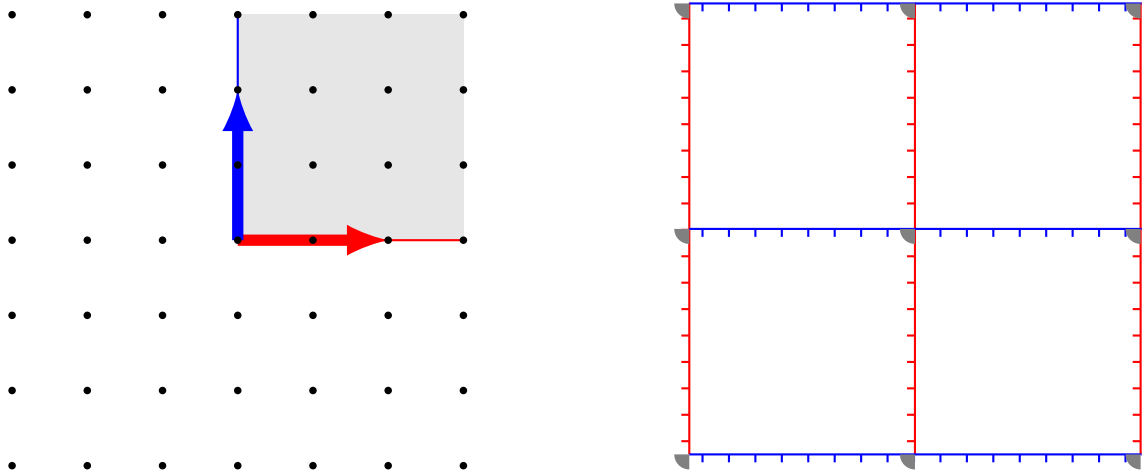

\begin{figure}[h]
\includegraphics[width=12cm]{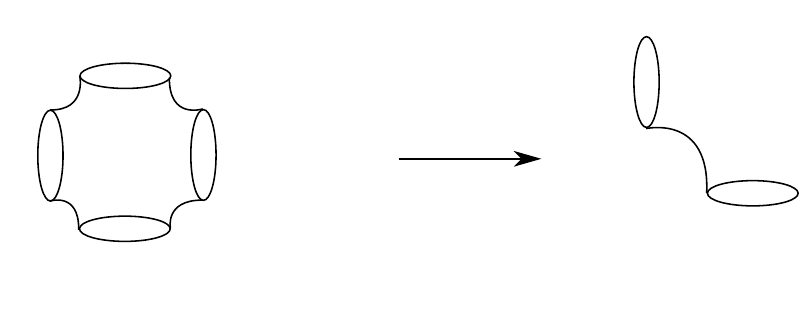}
\caption{The skeleton $\LL_{\Delta}$ of the pants $\tilcalP_\Delta$ associated to the simplex $\Delta$ with vertices $(0,0),(2,0),(0,2),$ and its 4:1 cover of the skeleton of the standard pants.}
\label{fig:stacky-cover}
\end{figure}
\end{example}

\section{Patchworking and skeleta} \label{sec:skeleta}

Fix a complex torus $\TT^\vee_\CC=T\TT^\vee$, along with a toric partial compactification 
$\bT_\Sigma$ arising from a (stacky) fan $\Sigma \subset M^\vee_\RR$.  We write 
$\Delta^\vee$ for the convex hull of the stacky primitives. 

According to \cite{HV}, the mirror to $\bT_\Sigma$ 
is the Landau-Ginzburg model associated to a function $W_\Sigma:\TT^\vee_\CC \to\CC$  
whose Newton polytope is $\Delta^\vee$.  In addition, the expected mirror
to $\partial \bT_\Sigma$ is a general fiber $F_W$ of $W_\Sigma$. 

In this section we will 
explain how $W_\Sigma$ determines a Liouville sector (i.e. prove Prop. \ref{prop:canembedfiber}) 
and  show that the relative skeleton of this sector is the FLTZ Lagrangian $\LL_\Sigma$.  

Let us briefly outline the ideas involved.  We will 
study the hyperplane $F_W$ through its \emph{amoeba} (\cite{GKZ}), the projection of $\dtmir$ to the tangent fiber: 
$$\calA:=\Log(\dtmir) \subset M_\RR.$$ 

The cones of $\Sigma$ give a triangulation of the polytope $\Delta^\vee$. 
We choose the Laurent polynomial $W_\Sigma$ so that its tropicalization $\Pi_\Sigma$ is a spine onto which $\calA$ retracts.
The complex $\Pi_\Sigma$ is a piecewise-affine locus dual to the triangulation of $\Delta^\vee$ by the cones of $\Sigma$. By 
assumption, this triangulation is {\em star-shaped} (all non-boundary simplices share a common vertex 0); 
the distinguished vertex corresponds to a distinuished component of the complement of the amoeba.
We denote the boundary of this region by $\partial^0\calA \subset \partial \calA$.  

Mikhalkin  \cite{M} shows how to isotope the hypersurface $F_W$ to another hypersurface $F_\Sigma$ whose amoeba 
is ``close" to the spine $\Pi_\Sigma.$
As we have recalled in Sections~\ref{sec:tailoring} and \ref{sec:skeleta2}, this isotopy was used by Nadler \cite{Nwms} to compute the skeleton of the ``$n$-dimensional  pants'',  {\em i.e.}, the zero locus of the 
polynomial $W_\Sigma = 1 + \sum_{i=1}^n z_i$.  

In our more general setting, 
Mikhalkin's isotopy ensures that the critical points of $\Log|_{F_\Sigma}$ -- and in fact the entire skeleton
 $\LL_\Sigma$ -- lie above the distinguished boundary component $\partial^0\calA$ of the amoeba.
The preimage of such a boundary component is precisely a contact type hypersurface.
Finally, to each pants in the decomposition of $F_\Sigma$ we apply the argument from \cite{Nwms} described in the previous section to obtain
the precise form of the skeleton.

\subsection{Pants decomposition of $F_\Sigma$}\label{sec:tropical-rec}
In order to construct $F_\Sigma$ and produce its skeleton, we will follow \cite{Nwms} in using Mikhalkin's theory of localized hypersurfaces, which we now recall. 

\subsubsection{Triangulation and dual complex}\label{subsec:triang}

Recall that we are assuming the fan $\Sigma$ is smooth and quasi-projective, or equivalently, that the subdivision $\mathcal{T}$ of $\Delta^\vee$ is a {\em regular} triangulation.
By definition, regularity of $\mathcal{T}$ means that $\mathcal{T}$ is the corner locus of a 
convex piecewise-linear function $\alpha:\Delta^\vee\cap M^\vee\to\RR$.
The \emph{Legendre transform} of $\alpha$ is the function

$$
\xymatrix{
L_\alpha:M_\RR\to \RR,&
m\mapsto\max_{n\in\Delta^\vee}(\langle m,n\rangle-\alpha(n)).
}
$$
\begin{defn}
	The \emph{dual complex} for the regular triangulation $\mathcal{T}$ is the polyhedral complex in $M_\RR$ obtained as the corner locus of the Legendre transform $L_\alpha$. We will denote the dual complex for $\mathcal{T}$ by $\Pi_\Sigma.$
\end{defn}

\begin{example}\label{ex:PLpants}Let $e_1,\ldots,e_n$ be a basis of $M^\vee,$ and let $\Delta_{std}^\vee$ be the polytope with vertices $0,e_1,\ldots,e_n.$ Then we can define a piecewise-linear function $\alpha$ on $\Delta_{std}^\vee$ by declaring $\alpha(0)=0,\alpha(e_i)=\alpha_i$ for some $\alpha_i>0$.  The resulting dual complex $\Pi_{std}$ is the corner locus of the function $(a_1,\ldots,a_n)\mapsto \max(0,a_1-\alpha_1,\ldots,a_n-\alpha_n)$; in other words, it is a translation by $(\alpha_1,\ldots,\alpha_n)$ of the tropical pants $\Pi_{n-1}$ defined in Section \ref{subsec:pants}.
\end{example}

The geometric significance of $\Pi_\Sigma$ is the following.  Recall that the \emph{amoeba} of a hypersurface in $(\CC^\times)^n$ is its image in $\RR^n$ under the map
$$
\xymatrix{
\Log:(\CC^\times)^n\to\RR^n,&\Log(z_1,\ldots,z_n)=(\log(|z_1|),\ldots,\log(|z_n|)).
}
$$

\begin{proposition} \cite{M}\label{prop:dom-trop}
Let $V$ denote the set of vertices in the triangulation $\calT,$ and let $H^t = \{f^t=0\},$ where
\[
f^t = \sum_{m\in V}t^{-\alpha(m)}z^m.
\]

For $t\gg0$, the complex $\Pi_\Sigma \cdot \log(t)$ will sit as a spine inside the amoeba of $\Log(H_t)$, 
and as $t\to \infty,$ the rescaled amoebae $\Log(H_t)/\log(t)$ converge (Gromov-Hausdorff) to $\Pi_\Sigma.$
\end{proposition} 
\begin{proof}
The basic idea is as follows.  Consider a face  $E$ of the dual complex $\Pi_\Sigma$, corresponding to a face $E^\vee$ of the triangulation $\calT$. 
Then the portion of the amoeba lying over the interior of $E$ is a region where the behavior of $f_t$ is dominated by those monomials in $f_t$ 
corresponding to vertices of $E^\vee$. See \cite[Section 6]{M} for details.
\end{proof}

One says that the complex $\Pi_\Sigma$ is the \emph{tropical hypersurface} associated to the Newton polytope $\Delta^\vee$ with regular triangulation $\mathcal{T}.$
We term $t$ the `tropicalization parameter'.

\begin{figure}[h]
\includegraphics[width=10cm]{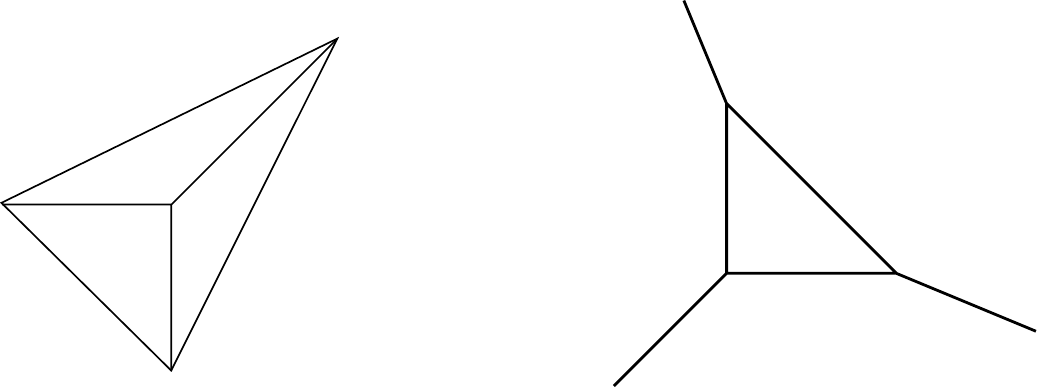}
	\caption{\label{fig:expi}A maximal subdivision of the polygon with vertices $(-1,0),(0,-1),(1,1)$, and its (unimodular) dual complex $\Pi.$}
\end{figure}

Since the triangulation $\mathcal{T}$ of $\Delta^\vee$ is star-shaped, we have: 

\begin{lemma}\label{lem:pi-str}
	Let $\Pi_\Sigma^0$ be the component of $M_\RR\setminus\Pi_\Sigma$ corresponding to the vertex 0 of the triangulation $\mathcal{T}$, and denote its boundary
	by $\partial^0\Pi_\Sigma.$ The polytope $\Pi^0_\Sigma$ is a (possibly unbounded) polytope with face poset anti-equivalent to the poset of nonzero cones in the fan $\Sigma.$
\end{lemma}
The polytope $\Pi_\Sigma^0$ will be bounded if and only if the toric variety $\bT_\Sigma$ is proper, in which case $\Pi_\Sigma^0$ will be the only bounded polytope in $M_\RR\setminus \Pi^0_\Sigma.$ 

\subsubsection{Tropical pants}
We required the subdivision $\mathcal{T}$ of $\Delta^\vee$ to be a \emph{triangulation}, which means that all of the faces in $\mathcal{T}$ are simplices. This allows us to divide up $\Pi_\Sigma$ into pieces we understand:
\begin{defn}
	The neighborhood in $\Pi_\Sigma$ of any vertex is a \emph{tropical pants}.
\end{defn}
These pants will be our basic building blocks in the construction to follow. This has two appealing features: the first is that the complex $\Pi_\Sigma$ is obtained by gluing these pants together.
Second, a $(k-1)$-face in $\Pi_\Sigma$ is the product of $\RR^k$ with a $(n-k-1)$-dimensional tropical pants. Hence the loci along which pants involved in the description of $\Pi_\Sigma$ are glued are products of the form $(\text{lower-dimensional pants})\times \RR^k$.

\subsubsection{Tailoring}
We now recall the construction of \cite{M} giving an isotopy from $F_W$ to some $F_\Sigma$ whose amoeba is closer
to the tropical hypersurface $\Pi_\Sigma$. 
In the case of the pants $\calP_{n-1},$ this isotopy was described in Proposition~\ref{prop:locfact} above.

It is straightforward to see what $F_\Sigma$ should be.  
Suppose two simplices $P_1,P_2$ in the triangulation $\calT$ share a common face $F$, so that their respective dual complexes $\Pi_{P_1},\Pi_{P_2}$ overlap in a common subcomplex $\Pi_F$, and
let $U$ be a neighborhood of the interior of $\Pi_F.$ Then the inductive structure of the tailored $P$-pants ensures that 
above $U$, the pants $\tcP_{P_1}$ and $\tcP_{P_2}$ agree: both are equal to the tailored leg $\tcP_{F}.$
Thus we may take the union of all these pants to define $F_\Sigma$. 

The isotopy can be glued similarly: 

\begin{lemma}[{\cite[Section 6.6]{M}, \cite[Propositions 4.2, 4.9]{A1}}]\label{lem:lhconstr}There is a Hamiltonian isotopy of
symplectic hypersurfaces $F_W \to F_\Sigma$ such that for each face $F$ in the tropical curve $\Pi_\Sigma\subset M_\RR$, corresponding to a polytope $P$ in the triangulation $\mathcal{T}$, there is a neighborhood $U_P\subset M_\RR$ such that 
	$\Log^{-1}(U_v)$ is equal to the intersection $\tilcalP_P^\ell$ with a large ball in $\TT_\CC^\vee.$
\end{lemma}
\begin{figure}[h]
\includegraphics[width=5cm]{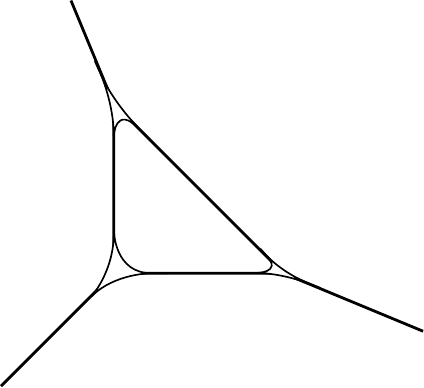}
	\caption{\label{fig:exloc}The amoeba of the localization of the hypersurface $xy+\frac{1}{x}+\frac{1}{y}=0.$}
\end{figure}
\begin{remark}\label{rem:isodef} The symplectic isotopy from \cite{A1} is defined as follows: for $t\geq0$ and $0\leq s\leq 1,$ write $H^{t,s} = \{f^{t,s} = 0\},$
  \begin{equation}\label{eq:tailoring-formula}
  f^{t,s}:=\sum_m t^{-\alpha(m)}(1-s\phi_m(\Log(z)))z^m,
\end{equation}
where the sum is taken over the vertices of the triangulation $\calT$ of $\Delta^\vee,$ and $\phi_m\in C^\infty (\RR^n)$ is a certain function which is 1 in a neighborhood of the component of $M_\RR\setminus \Pi_\Sigma$ corresponding to $m,$ and 0 away from that region; as in Section~\ref{subsec:triang}, taking the tropicalization parameter $t$ large ensures that $\Log(H^{t,s})$ contains $\log(t)\cdot \Pi_\Sigma$ as a spine. Taking the ``tailoring parameter'' $s$ from 0 to 1 deforms the hypersurface $\{\sum_m t^{-\alpha(m)}z^m\}$ by forcing that, on each region of the amoeba $t$, any term which does not dominate the behavior of $f^{t,0}$ in that region (as described in the proof of Proposition~\ref{prop:dom-trop}) does not contribute at all. 
\end{remark}
\begin{remark}
  As in our definition of the standard pants, our convention in this paper will differ from that in Equation~\eqref{eq:tailoring-formula} by our choice to take the sign of the constant coefficient of $f^{t,s}$ to be negative rather than positive. This ensures that the real positive points of $H^{t,s}$ lie over the boundary of the central component of the amoeba complement.
\end{remark}

\subsection{The skeleton of $F_\Sigma$}\label{subsec:main}
As in Section~\ref{sec:skeleta2}, by choosing an inner-product on $M$, we obtain an isomorphism $\TT_\CC^\vee=T\TT^\vee\cong T^*\TT^\vee,$ and we restrict the symplectic form $\omega,$ and its primitive $\lambda$ from this space to $F_\Sigma$. 
We will use the pants decomposition of $F_\Sigma$ to (observe that it is a Liouville manifold and) compute its skeleton, which we denote by $\Lambda_\Sigma.$ 

However, in order to avoid performing any calculations beyond those described so far, we must adopt
a certain technical hypothesis on the fan $\Sigma$.  As remarked in the introduction, this hypothesis can be removed;
see \cite{Z} for details.   

\begin{defn}
	A polytope $P\subset M_\RR^\vee$ is called \emph{perfectly centered} if for each nonempty face $F\subset P,$ the normal cone of $F$ (transported to $M_\RR^\vee$ by the inner product $M_\RR^\vee\cong M_\RR$) has nonempty intersection with the relative interior of $F$.
\end{defn}

As in the proof of Lemma~\ref{lem:regtri}, we write $\alpha:\Delta^\vee\to\RR$ for a function inducing the regular triangulation of $\Delta^\vee$ 
defined by $\Sigma.$ The complex $\Pi_\Sigma$ depends on our choice of $\alpha$.

\begin{defn}
	We will say that a fan $\Sigma$ is \emph{PC} if there exists some $\alpha$ as above for which the polytope $\Pi_\Sigma^0$ is perfectly centered.
\end{defn}
Assume now that the fan $\Sigma$ is PC.

\begin{remark}
	So far, no fan is known to us not to be PC; nor, however, do we know any compelling reason why all fans should be PC.
\end{remark}

We will denote the amoeba of $F_\Sigma$ by 
\[\tcA_\Sigma:=\Log(F_\Sigma).\]
Recall that we write $\Pi^0_\Sigma$ for the component of $M_\RR\setminus \Pi_\Sigma$ dual to the unique 0-dimensional cone in $\Sigma$ and $\partial^0\Pi_\Sigma$ for its boundary. Write $\partial^0\tcA_\Sigma$ for the corresponding boundary component of the amoeba, and

Recall that we write $-\LL_\Sigma=\bigcup_{0\neq \sigma\in \Sigma}\sigma^\perp\times \sigma$ for the (negative) FLTZ skeleton.

\begin{figure}
	\centering
	\includegraphics[width=5cm]{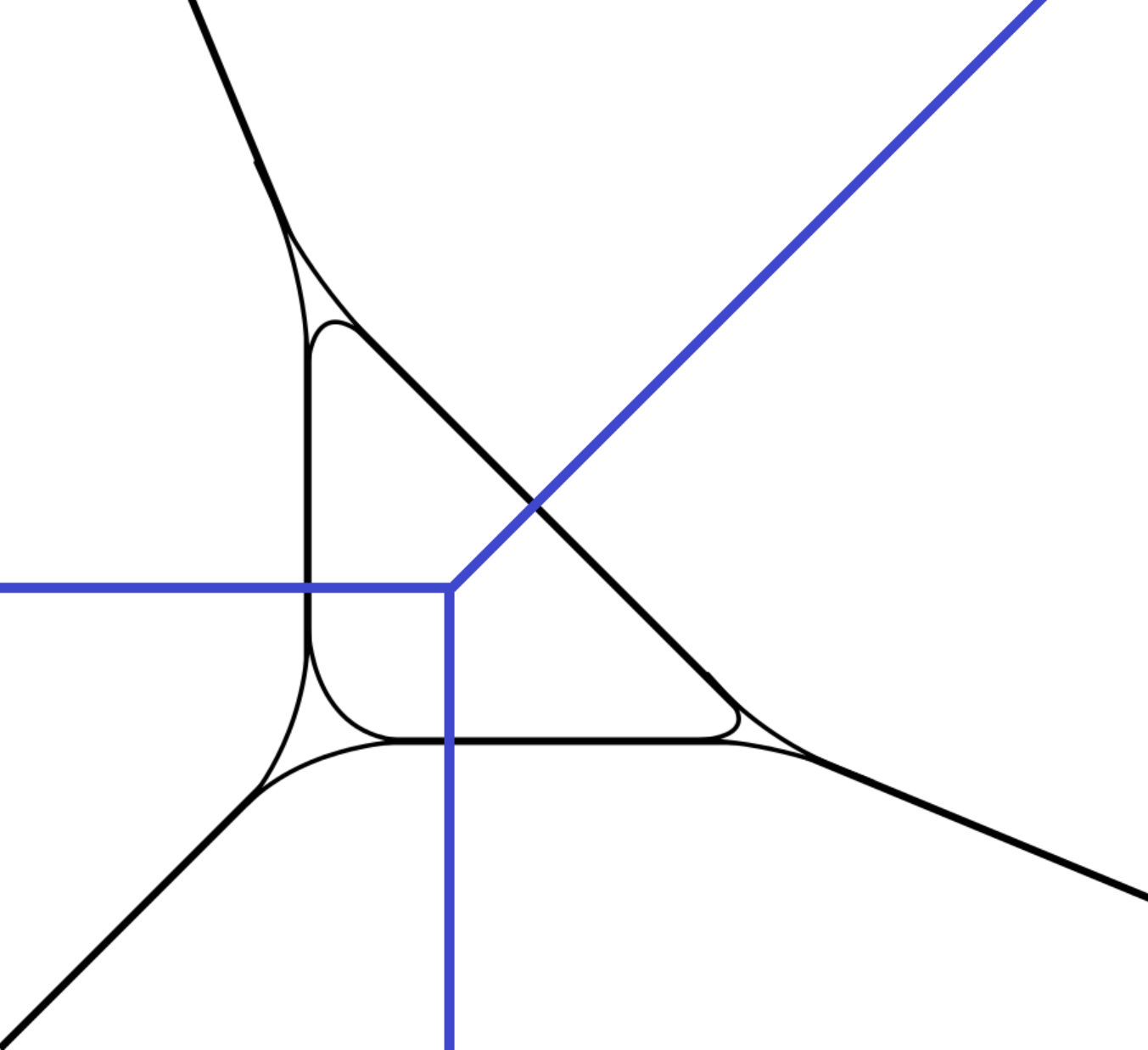}
	\caption{The fan $\Sigma$ for $\PP^2$ superimposed on the amoeba $\calA_\Sigma.$} 
	\label{fig:complete-example}
\end{figure}

\begin{theorem}\label{thm:hypskel}
The skeleton $\Lambda_\Sigma$ of $\ttmir$ can be written as the intersection
$$
\xymatrix{
\Lambda_\Sigma=C\cap(-\LL_\Sigma).
}
$$
\end{theorem}
\begin{proof}
	From our hypothesis that the fan $\Sigma$ is PC, we may assume that the polytope $\Pi_\Sigma^0$ is perfectly centered,
	so that each nonzero cone $\sigma$ in $\Sigma$ intersects its dual face in $\Pi_\Sigma^0$, as in Figure~\ref{fig:complete-example}.
	This allows us to define an open cover of $F_\Sigma$ as follows: for each top-dimensional cone $\sigma$ in $\Sigma$, let $V_\sigma$ be a neighborhood of the cone $\sigma,$ thought of as in $M_\RR$. 

	Let $U_\sigma=\Log^{-1}(V_\sigma)\cap F_\Sigma$ be the lift of $V_\sigma$ to an open subset of $F_\Sigma.$ Then $U_\sigma$ is an open subset in a pants $\tilcalP_\sigma^\ell$. By construction, the image of $U_\sigma$ in $\tilcalP_\sigma^\ell$ contains the whole skeleton $\Lambda_\sigma$ of $\tilcalP_\sigma^\ell.$
On the other hand, every zero of $\lambda|_{F_\Sigma}$ is contained in some $U_\tau,$ as is its stable manifold.  We conclude the skeleton $\Lambda_\Sigma$ is equal to the union of the skeleta $\Lambda_\tau.$
\end{proof}

Note $C$ is transverse to the Liouville flow on $\TT^\vee$, hence contact.  In addition:  

\begin{lemma}\label{lem:nottangent} In a neighborhood of the skeleton $\Lambda_\Sigma,$ the hypersurface $F_\Sigma$ is nowhere tangent to the ambient Liouville vector field of the Weinstein manifold $\TT_\CC^\vee.$
\end{lemma}
\begin{proof}
  The pants cover of $F_\Sigma$ allows us to reduce to the case where $F_\Sigma = \tcP_P^\ell$ is a $P$-pants. Now we can proceed by the induction used in the proof of Lemma~\ref{lem:nonstd}. 

  In the base case, $P=\Conv(0,v)$ is 1-dimensional, and  $\tcP_P^\ell=\{z_1^{v_1}\cdots z_n^{v_n}=e^\ell\}$ is a copy of $(\CC^\times)^{n-1}$ projecting by Log to its tropical hypersurface $\Pi_\Sigma$. The Liouville vector field on $\TT_\CC^\vee=(\CC^\times)^n$, under the Log projection, points directly outward from $\Pi_\Sigma.$

  Now suppose $P=\Conv(0,v_1,\ldots,v_n)$ is an $n$-simplex. We know the result on the legs of $\tcP_P^\ell$ by induction, so we only need to prove it in a neighborhood of the big simplex $S_+$ in the skeleton, which is the preimage under the cover $\tcP_P^\ell \to \tcP_{n-1}^\ell$ of the positive real points of $\tcP_{n-1}^\ell.$ But the Liouville vector field on $\TT_\CC^\vee=(\CC^\times)^n$ in coordinates $z_j=e^{\xi_j+i\theta_j}$ is $\sum_j \xi_j\partial_{\xi_j},$ so if $\tcP_P^\ell$ had any tangent vectors along $S_+$ in the direction of the Liouville flow, this would imply that $\tcP_{n-1}^\ell$ had positive real points which do not project to the boundary of the amoeba, which is false.
\end{proof}

\begin{cor}
There is a Liouville domain $D \subset \TT^\vee$ completing to $\TT^\vee$ and a Liouville domain 
$F \subset F_\Sigma$ completing to $F_\Sigma$, such that $F \subset \partial D$ and 
the FLTZ Lagrangian $-\LL_\Sigma$ is a relative skeleton for the pair $(D, F)$. 
\end{cor}
\begin{proof}
The point is that we may deform $C$ transversely to the Liouville flow in such a way as to 
cause it to contain some neighborhood $U_\Sigma \subset F_\Sigma$ of $\Lambda_\Sigma$.  

Indeed, the Liouville flow gives an identification $T^* \TT^\vee \setminus \TT^\vee \cong C \times \R$, with
$C$ included as $C \times 0$.  By 
Lemma \ref{lem:nottangent}, for some closed manifold $V_\Sigma \subset F_\Sigma$ neighborhood of $\Lambda_\Sigma$, 
the corresponding projection $V_\Sigma \to C$ is an embedding.  Its image is some codimension one 
smooth hypersurface (with boundary) of $C$, over which we may write $V_\Sigma$ as the graph 
of a smooth function.  Extend this function arbitrarily to all of $C$. The graph of the result will be the boundary
of our desired $D$. 

We thank John Pardon for this method of constructing Liouville pairs.
\end{proof}

\begin{example}
	Let $\Delta^\vee$ be the polytope with vertices $(1,1),(0,-1),(-1,0),$ as in Figure~\ref{fig:expi} and Figure~\ref{fig:exloc}. In Figure~\ref{fig:complete-example}, the fan $\Sigma$ is drawn superimposed on the amoeba $\calA_\Sigma.$ A neighborhood of each top-dimensional cone in $\Sigma$ is a pair of pants which contributes to $\LL$ a pair of circles attached by an interval. The circles live over the points where the rays of $\Sigma$ intersect $\Pi,$ and the intervals lie over the boundary of the bounded region in the center of the amoeba.
\end{example}

\section{Microlocalizing Bondal's correspondence}  \label{sec:microccc}

Recall we denote by $\TT$ a torus with respective character and
cocharacter lattices $M$ and $M^\vee$.   Fix a (stacky) fan $\Sigma \subset M^\vee_\RR$ 
and the corresponding toric partial compactification $\TT_\CC \subset \bT_\Sigma$. 

Bondal \cite{B} described a fully faithful embedding
of the category of coherent sheaves on $\bT_\Sigma$ into the 
category of constructible sheaves on the real torus $\TT^\vee_\RR := M \otimes \RR/\ZZ$. 
This was developed further in \cite{FLTZ2, FLTZ3, Tr}; in particular, the constructible
sheaves in question were observed to have microsupport contained in $\LL_\Sigma$ and conjectured
to generate the category of such sheaves.  This conjecture was established in \cite{Ku}.

We use this equivalence to prove a similarly-flavored equivalence ``at infinity",
\emph{i.e.}, an equivalence between the category of coherent sheaves on the 
toric boundary and the category of wrapped microlocal sheaves away from the zero section.

\bigskip
\emph{Categories and conventions} ---  We work with dg categories over a fixed ground ring $k$.  This theory can
be set up either directly 
\cite{Kel-deriving, Kel-survey, Dr} or by specializing the theory of stable $(\infty, 1)$-categories of \cite{Lur-topos, Lur-algebra} as in  \cite[I.1.10]{GR}.  

The microlocal sheaf theory of \cite{KS} was originally developed in 
the setting of the bounded derived category.  It is essential for our work here to work with the 
dg category of unbounded complexes.  It is well known to experts that it is straightforward
to set up the sheaf theory in this setting (see e.g. \cite[Sec. 2.2]{Ncs} or \cite[Sec. 4.1]{GPS3}) 
and that, with the use of \cite{Spa, RS} to deal with some issues around unbounded complexes, 
all constructions of \cite{KS} may be translated to this setting. 
 
For a manifold $M$, we write $Sh(M)$ for the unbounded dg derived category of sheaves of $k$-modules on $M$.  
We impose no restrictions on the stalks; \emph{i.e.}, we write $Sh$ for what in \cite{Nwms} is called $Sh^\diamond$ (and similarly
for the later $\mu sh$).

For a conical subset $Z \subset T^*M$, we write $Sh_Z(M)$ for the full subcategory of $Sh(M)$ consisting of those sheaves
with microsupport in $Z$.  When $Z$ is subanalytic Lagrangian, then this subcategory is compactly generated,
and we write $Sh_Z(M)^c$ for the subcategory of compact objects.  This subcategory is generally larger than the category of sheaves
with perfect stalks in $Sh_Z(M)$; for instance, when $Z = \emptyset$ it contains the tautological (derived) local system with fiber 
$C_* (\Omega M)$.  The idea to use compact objects in the unbounded category to model the {\em wrapped} Fukaya category stopped at $Z$
is due to Nadler \cite{Nwms}; that it works is now a theorem \cite{GPS3}.  The reader is referred to these articles for further discussions
of this category.

For $X$ an algebraic variety (or stack), we write $\QCoh(X)$ for the dg derived category of quasi-coherent sheaves on $X$ in
the sense of \cite{GR}; as observed there, the bounded subcategory agrees with the usual usage of this term.  
It is useful to remember that perfect complexes (bounded complexes of projectives)
are precisely the compact objects in $\QCoh(X)$, which can be recovered from $\Perf(X)$ by ind-completion. 
Similarly, we will write $\IndCoh(X)$ for the Ind-completion of the category $\Coh(X)$ of coherent sheaves on $X$ (\cite{GR}). We can recover
the category $\Coh(X)$ by passing to compact objects.

To simplify notation, we write as if $\Sigma$ is an ordinary (non-stacky) fan. 
To arrive at the corresponding statements in the stacky case, one need merely remember 
the data of a finite abelian group $\Gamma_\sigma$ for each cone in $\sigma$, and correspondingly 
replace the sets
$\{A(\sigma)\}_{\sigma\in\Sigma},$ $\{B(\sigma)\}_{\sigma\in\Sigma}$ with sets $\{A(\sigma,\chi)\}_{\sigma\in\Sigma,\chi\in \Gamma_\sigma}$, $\{B(\sigma,\chi)\}_{\sigma\in\Sigma,\chi\in \Gamma_\sigma},$ where the added $\chi$ denotes translation in $\TT^\vee$ and twists by a character, respectively. See \cite[Section 5]{FLTZ3} for details.

\subsection{Bondal's coherent-constructible correspondence}

$ $ \vspace{2mm}

For a cone $\sigma \subset M^\vee$, we write 
$B(\sigma)$ for the structure sheaf on $\Spec(k[\sigma^\vee])$, or its pushforward
to any toric variety whose fan contains the cone $\sigma$.  On the other hand,
we write $A(\sigma)$ for the constructible sheaf on $M \otimes \RR/\ZZ$ obtained
by taking the $!$-pushforward of the dualizing (constructible) sheaf on the interior
of $\sigma^\vee$. One then makes the following

\vspace{2mm}

{\bf Basic calculation} (\cite{B, FLTZ2, Tr})\textbf{:}
Let $\bT_\Sigma$ be a toric variety with fan $\Sigma$, with dense torus $\TT_\CC.$  Let $\sigma, \tau\in\Sigma$ be cones.  Then there
are canonical isomorphisms
$$H^*\Hom(A(\sigma), A(\tau))  \cong k[\tau^\vee] \cong H^* \Hom(B(\sigma), B(\tau))\qquad \qquad \mbox{if}\,\, \sigma \supset \tau,$$
and all other Homs between such objects vanish.  This is moreover compatible with the evident composition structure.

\vspace{2mm}
We denote full dg subcategories generated by the $A(\sigma)$ and $B(\sigma)$ by
\begin{eqnarray*}
 A_\Sigma & := & \{A(\sigma)\, |\, \sigma \in \Sigma\} \subset Sh(\TT^\vee), \\
 B_\Sigma & := & \{B(\sigma)\, |\, \sigma \in \Sigma\} \subset QCoh(\bT_\Sigma).
\end{eqnarray*}

While the calculation above might seem to imply only the equivalence $H^0(A_\Sigma) \cong H^0(B_\Sigma)$ of \emph{triangulated} categories,
we recall the following useful fact: 

\begin{lemma} \label{lem:freecoherence}
Let $\mathcal{C}_i$ be a collection of dg categories, each of which has all morphisms concentrated in cohomological 
degree zero.  Then any diagram valued in the $H^0(\mathcal{C}_i)$ lifts canonically to a homotopy coherent diagram
in the corresponding $\mathcal{C}_i$. 
\end{lemma}
\begin{proof}
The hypothesis on $\calC_i$ implies that the natural maps 
$$
\xymatrix{
H^0(\mathcal{C}) &\ar[l] \tau_{\le 0} \mathcal{C} \ar[r]& \mathcal{C}
}
$$
are quasi-isomorphisms.  Thus any diagram among the $H^0(\calC_i)$ can be lifted to a diagram among the $\calC_i$ by composing
with this pair of quasi-isomorphisms. 
\end{proof}

As the category of quasicoherent sheaves on a toric variety is generated by the structure sheaves of the affine toric charts, 
the restriction to the subcategory $B_\Sigma$ is really no restriction: the morphism $QCoh(\bT_\Sigma) \to \Mod-B_\Sigma$ 
is an isomorphism.  

On the other side, the objects of $A_\Sigma$ all satisfy the microsupport estimate 

$$ss(A(\tau)) \subset \bigcup_{\sigma \subset \tau} \sigma^\perp \times (-\sigma) \subset \TT^\vee \times M^\vee_\RR = T^* \TT^\vee.$$
In particular, writing
$$\LL_\Sigma := \bigcup_{\sigma \in \Sigma} \sigma^\perp \times (-\sigma),$$
we have that $A_\sigma \in Sh_{\LL_\Sigma}(\TT^\vee) $ for all $\sigma \in \Sigma$. 
As conjectured by \cite{FLTZ, Tr}, and proven by Kuwagaki \cite{Ku}, these objects generate this category:  

\begin{theorem} \cite{Ku}  \label{thm:Kuwagaki}
When $\bT_\Sigma$ is a smooth orbifold, the morphism $Sh_{\LL_\Sigma}(\TT^\vee) \to \Mod-A_\Sigma$ is an isomorphism.
\end{theorem} 

\begin{remark} \label{rem: Kuwagaki}
In fact what Kuwakagi proves is that for any, not necessarily smooth, $\bT_\Sigma$ there is an isomorphism 
$Sh_{\LL_\Sigma}(\TT^\vee) \cong IndCoh(\bT_\Sigma)$.  The above statement follows because 
in the smooth case, $IndCoh = QCoh$, which as we mentioned above is generated by the $B_\sigma$.
We use the above formulation rather than Kuwagaki's more general result because we will later make calculations
with the $A_\sigma$ and $B_\sigma$ directly, and   
we restrict ourselves to the smooth case to avoid e.g. worrying about how to lift the $B_\sigma$ to $IndCoh$. 

This causes no loss of generality, since it is anyway only 
in the smooth case that we have been able to identify $\partial \LL_\Sigma$ as 
a relative skeleton. 
\end{remark}

\subsection{Restriction is mirror to microlocalization}
Let $\bT$ be a toric variety, $\sigma$ a cone of the fan $\Sigma(\bT)$, 
and $i_\sigma: \overline{O(\sigma)} \to \bT$ the inclusion of the orbit closure corresponding 
to the cone $\sigma$.   As the orbit closure is itself a toric variety, one can ask
what functor of constructible sheaf categories corresponds under Bondal's correspondence
to the pullback $i_\sigma^*.$
We will see that the answer is a sort of microlocalization functor.

\subsubsection{Restriction to orbit closures}

Recall that the orbit closure $\overline{O(\sigma)}$ carries the structure of a toric variety, with 
associated cocharacter lattice $M^\vee/\ZZ\sigma$.  For $\tau$ a cone containing $\sigma$,
we write $\tau/\sigma$ for the image of $\tau$ in $M^\vee/ \ZZ \sigma$. 
The map $\tau \to \tau/\sigma$ gives a bijection between cones containing $\sigma$
and cones in the fan of $\Sigma(\overline{O(\sigma)})$.   We will therefore write
$\Sigma / \sigma := \Sigma(\overline{O(\sigma)})$. 

Let us recall that

$$\bT_\tau = \coprod_{\tau \supset \eta } O(\eta)$$
$$\overline{O(\sigma)} = \coprod_{\eta \supset \sigma} O(\eta),$$
and therefore the intersection of the orbit closure $\overline{O(\sigma)}$ with the affine piece $\bT_\tau$ decomposes as
$$\bT_\tau \cap \overline{O(\sigma)} = \coprod_{\tau \supset \eta \supset \sigma} O(\eta) = \begin{cases} 
\overline{O(\sigma)}_{\tau/\sigma} \qquad \qquad \tau \supset \sigma, \\ \emptyset \qquad \qquad \qquad
\mbox{otherwise}. \end{cases}$$

For $\tau \supset \sigma,$ there is a natural identification
$(\tau/\sigma)^\vee \cong \tau^\vee \cap \sigma^\perp \subset \tau^\vee$.  The corresponding map 
$k[(\tau/\sigma)^\vee] \hookrightarrow k[\tau^\vee]$ has a unique $M$-graded left-inverse
$k[\tau^\vee] \to k[(\tau/\sigma)^\vee]$, which gives the affine inclusion 
$\overline{O(\sigma)}_{\tau/\sigma} \hookrightarrow \bT_\tau$.  We conclude:

\begin{lemma}
We have canonical isomorphisms
$$i_\sigma^* B(\tau) = \begin{cases} B(\tau/\sigma) \qquad \qquad  \tau \supset \sigma, \\ 0 \qquad \qquad \qquad \mbox{otherwise}. \end{cases} $$

The source or target of the induced map  $i_\sigma^*: H^*\Hom(B(\tau'), B(\tau)) \to H^*\Hom(B(\tau' / \sigma), B(\tau /\sigma))$ vanishes unless
$\tau' \supset \tau \supset \sigma$, and in this case is canonically identified with the map
$k[{\tau}^\vee] \to k[ (\tau/\sigma)^\vee]$. 
\end{lemma}

\subsubsection{Microlocalization}
Our description of the mirror to the restriction functor $i_\sigma^*$ will be given in terms of Sato's microlocalization. 
We now briefly review this notion; for details
see \cite[Chap. 4]{KS}. 

Microlocalization is built from Verdier specialization, 
and the Fourier-Sato transform.  The Verdier specialization along a submanifold $X \subset Y$
carries sheaves on $Y$ to conic sheaves on $T_X Y$, by pushing forward along a deformation to the normal cone.  The Fourier-Sato
transformation carries conic sheaves on a bundle to conic sheaves on its dual, by convolution with the kernel given by the constant
sheaf on the locus $\{(x, x^*)\, |\, x^*(x) \le 0\}$.  Sato's microlocalization is the composition of these, and carries sheaves on $Y$ to
conic sheaves on $T_X^* Y$; we denote it by $\mu_X$. 

As usual, write
$$\LL_{\Sigma} = \bigcup_{\sigma \in \Sigma} \sigma^\perp \times (-\sigma) \subset T^* \TT^\vee$$ 
for the \cite{FLTZ} skeleton mirror to $\bT_\Sigma$.

For the orbit closure
$\overline{O(\sigma)}$, we denote the corresponding torus $\TT_\CC(\sigma) := \TT_\CC / (\ZZ\sigma \otimes \G_m)$, and the corresponding skeleton
$\LL_{\Sigma/\sigma} \subset T^*\TT(\sigma)^{\vee}$.  
Note the canonical identification $\TT(\sigma)^{\vee} \cong \sigma^\perp$.

We compute the following microlocalization: 

\begin{lemma}  \label{lem:CCCmicro}
Let $\pi: \sigma^\perp \times (-\sigma^{\circ}) \to \sigma^\perp \cong \TT(\sigma)^{\vee}$ be the projection.  
Then the morphism

\begin{eqnarray*}
m_\sigma: Sh(\TT^{\vee}) & \to & Sh(\sigma^\perp) \\
F & \mapsto & \pi_* ((\mu_{\sigma^\perp} F)|_{\sigma^\perp \times (-\sigma^\circ)}).
\end{eqnarray*}

respects FLTZ skeleta, i.e. restricts to $m_\sigma: Sh_{\LL_\Sigma} (\TT^{\vee}) \to Sh_{\LL_{\Sigma/\sigma}}(\sigma^\perp)$.  Moreover, there are canonical isomorphisms
$$m_\sigma A(\tau) = \begin{cases} A(\tau/\sigma) \qquad \qquad  \tau \supset \sigma \\ 0 \qquad \qquad \qquad \mbox{otherwise}. \end{cases} $$

The source or target of the induced map $\mu_\sigma: H^*\Hom(A(\tau'), A(\tau)) \to H^*\Hom(A(\tau' / \sigma), A(\tau /\sigma))$ vanishes unless
$\tau' \supset \tau \supset \sigma$; in this case the map is canonically identified with 
$k[{\tau}^\vee] \to k[ (\tau/\sigma)^\vee]$.
\end{lemma} 
\begin{proof}
As we will we see the sheaves in question are constant along the fibers of $\pi$, which are contractible, the pushforward $\pi_*$ does essentially nothing,
and we subsequently omit it from the notation. 

As $\LL_\Sigma$ is the union of the microsupports of the $A(\sigma)$, our argument showing that 
$m_\sigma A(\sigma) = A(\sigma/\tau)$ will also show that $m_\sigma(Sh_{\LL_\Sigma} (\TT^{\vee})) 
\subset Sh_{\LL_{\Sigma/\sigma}}(\sigma^\perp)$.  One could also directly use the formula
\cite[Thm. 6.4.1]{KS} showing that the microsupport of the microlocalization is the specialization of
the microsupport to the normal bundle of the conormal bundle. 

The vanishing when $\sigma\not\subset\tau$ follows immediately from the microsupport estimate (\cite[Prop. 5.1]{FLTZ})

$$ss(A(\tau)) \subset \bigcup_{\sigma \subset \tau} \sigma^\perp \times (-\sigma) \subset \TT^\vee \times M^\vee_\RR = T^* \TT^\vee.$$

Now consider $A(\tau)$ with $\sigma \subset \tau$.  The specialization of $A(\tau)$ along $\sigma^\perp$ 
can be understood as follows.  Choose a splitting $\TT^\vee = \sigma^\perp \times \TT'$, where $\TT' = \Hom(\ZZ \sigma, \RR/\ZZ)$.   
Let $\TT'_\epsilon$ be an epsilon ball around the origin of $\TT'$.  Then the Verdier specialization along $\sigma^\perp$ 
can be visualized as first restricting to $\sigma^\perp \times \TT'_\epsilon$,
and then rescaling the $\TT'_\epsilon$ factor to be very large, in the limit as $\epsilon \to 0$.  In this limit, the 
$\TT'_\epsilon$ factor can be identified with $\Hom(\ZZ \sigma, \RR)$.  

Restricting to $\sigma^\perp \times \TT'_\epsilon$ breaks $A(\tau)$ into a direct sum of $\NN^k$ pieces, 
where the $\NN^k$ grading counts how many times the cone has wrapped around $(S^1)^k$. 
Let us call the result $A'(\tau)$. 

First we study the grading zero component, $A'(\tau)_0$.  The rescaling limit carries $A'(\tau)_0$ to 
$A'(\tau)_0|_{\sigma^\perp} \boxtimes A_\epsilon(\sigma)$, where $A_\epsilon(\sigma)$ is the costandard sheaf on 
the dual cone to $\sigma$ inside $\Hom(\ZZ \sigma, \RR)$.  
The Fourier transform (which happens only in the second factor) of $A_\epsilon(\sigma)$ returns the standard sheaf on
$-\sigma$, which restricts to the constant sheaf on $-\sigma^\circ$.  On the other hand, $A'(\tau)_0|_{\sigma^\perp}$ 
is readily seen to be $A(\tau/\sigma)$. 

For the remaining components, note that since each has already wrapped around at least once in some direction, 
they are invariant along the line spanned by 
some extremal ray of the dual cone to $\sigma$ inside $\Hom(\ZZ \sigma, \RR/\ZZ)$.  It follows that their Fourier transform 
is supported on the face of $\sigma$ annihilated by that ray; hence the restriction of such a component to 
$-\sigma^\circ$ is zero. 

Finally, for the Homs, the above statement follows from the fact that $\NN^k$ grading coming from counting 
wrapping is identified with the natural gradings on $k[\tau^\vee]$, $k[(\tau/\sigma)^\vee]$, etc., as in \cite[Prop 2.3]{Tr}.
\end{proof}

In words: Bondal's correspondence intertwines the pullback $i_\sigma^*$ with the microlocalization
$m_\sigma$, at least as far as $A_\Sigma$ and $B_\Sigma$ are concerned.  By Theorem \ref{thm:Kuwagaki}, 
(and noting again Lemma \ref{lem:freecoherence}),
this can be extended to the larger categories.

\begin{remark} 
In \cite{FLTZ2, Tr} a different functoriality statement is established, which however does 
not apply to the case of an inclusion of a toric divisor.  Their result concerns morphisms
which, on the A-side, can be described in terms of just sheaves on the base manifold, rather than 
in terms of microlocalization.  
\end{remark}

\subsection{Microlocal sheaves}

\subsubsection{The Kashiwara-Schapira stack}
\label{subsub:msh}

Let $M$ be a manifold.  Using the tools of \cite{KS}, one can construct a 
sheaf of categories on $T^*M$, the
\emph{Kashiwara-Schapira stack}, whose global sections 
recover the usual category of sheaves on $M$.  To define it, one begins with the presheaf of categories $\mu sh^{pre},$ whose sections in a small ball $U$ are the quotient category
$$\mu sh^{pre}(U) = Sh(M) / Sh_{T^* M \setminus U} (M).$$
For a conical subset $\LL \subset T^*M$, there is a presheaf of full subcategories 
$\mu sh_\LL^{pre}$ on objects whose microsupport near $\LL$ is contained in $\LL$.

The Kashiwara-Schapira stack is the sheafification of this presheaf of categories; {\emph{i.e.}, it is obtained by
replacing sections by their limits over certain open covers.   While this sheaf of categories is not discussed 
explicitly in \cite{KS}, one does find there discussed its stalks (the $D^b(X; p)$ of \cite[Chap. 6]{KS}) and 
the Hom sheaves between global objects (under the name $\mu hom$).   
The actual sheaf of categories is discussed in some detail in \cite{Gui, Gui2, Nwms, NS}.\footnote{
The discussion in \cite[Chap. 10]{Gui2} gives many details, including an explanation of how the results of \cite{KS} may be translated into the assertions that $D^b(X; p)$ is the stalk of $\mu sh$ and that $\mu hom$ is the 
sheaf of homs.  While the official language of \cite{Gui2}  
is that of triangulated categories rather than DG categories -- an unfortunate choice insofar as triangulated structures do not glue well whereas DG categories do -- 
the constructions of \cite{Gui2} are all compatible with DG enhancement, and the proofs go through unchanged
once one has set up the basic sheaf theory in the DG setting.  These categorical aspects are 
discussed explicitly in \cite{Nwms, NS}.}  

To be precise,
let us specify in which $(\infty, 1)$-category of dg categories
these limits should be understood.  We use the following notation: 

We write $dg$ to mean the category whose objects are small stable (aka pre-triangulated) dg categories, and whose morphisms are exact functors.  
We write $DG$ for the category whose objects are cocomplete stable dg categories, and whose morphisms are exact functors.  
There are various not full subcategories of $DG$ characterized by what sort of adjoints the morphisms are.  We indicate
by ${}^* DG$ the category in which all morphisms are left adjoints; by ${}^* {}^*  DG$ the category in which all morphisms are left
adjoints of left adjoints; and so on.

Taking adjoints gives equivalences of categories switching the restrictions on adjoints; for instance, ${}^* DG \cong (DG {}^*)^{op}$, and so on.  This turns out to be very useful:
as described in \cite{Ga}, we can turn colimits into limits. 
Taking ind-completion and then adjoints gives an equivalence $dg \hookrightarrow {}^* {}^* DG \cong ({}^* DG {}^*)^{op}$; 
with the image being the compactly generated categories.  Thus
a colimit in $dg$ becomes a limit in ${}^* DG {}^*$, which we can compute in $DG^*$.
Taking adjoints again and passing to compact objects gives the originally desired colimit. 

We sheafify $\mu sh_\LL^{pre}$ in the $\infty$-category of $\infty$-categories.  The restriction
maps of $\mu sh_\LL^{pre}$ are continuous and cocontinuous; as such we could have
viewed the presheaf to be valued in any of $DG^*$, ${}^* DG$, ${}^* DG^*$, and sheafified there. 
But at least when $\LL$ is subanalytic Lagrangian (the only case of concern here),
in fact it does not matter where we sheafify, since the sections of the presheaf stabilize in a sequence of contractible open
neighborhoods around any given point $p \in \LL$.  In particular our $\mu sh_\LL^{pre}$ can be regarded
as a sheaf valued in ${}^* DG^*$. 

This sheaf becomes a cosheaf 
via the equivalence $({}^* DG {}^*)^{op} \cong  {}^* {}^* DG$.  Again because $\LL$ is subanalytic Lagrangian,
all sections are compactly generated; we may
pass to compact objects in the cosheaf to obtain a cosheaf valued in $dg$. 

\subsubsection{Microlocal restriction}

We now give some lemmas about how to compute the restriction of $\mu sh$. 
Let $X$ be a manifold and $M\subset X$ a submanifold.  We write $T^*_M X \subset T^*X$ 
for the conormal bundle to $X$.  Recall that the Sato microlocalization is a functor

$$\mu_M: Sh(X) \to Sh(T^*_M X)^{\RR^+}$$ 

Here the $\RR^+$ is to indicate that the sheaves are conic, \emph{i.e.}, constant along the cotangent directions.

Note that, locally near $T^*_M X$, the ambient cotangent bundle $T^*X$ is also the cotangent bundle of 
$T^*_M X$.  Thus it is natural to expect an expression for $SS(\mu_M F) \subset T^*(T^*_M X)$ in terms
of $SS(F) \subset T^*X$.  These cannot be equal in general, since $SS(\mu_M F)$ must be conic
in both the cotangent and in the cotangent-to-cotangent directions in $T^*(T^*_M X)$, whereas 
$SS(F) \subset T^*X$ will only be conic in the $T^*X$ cotangent directions.  This, however, 
is the ``only'' difference: 

\begin{theorem} \label{microsupport of microlocalization} \cite[Thm. 6.4.1]{KS}
$SS(\mu_M F) \subset T^*(T^*_M X)$ is obtained by specializing $SS(F)$ to the normal cone to 
$T^*_M X$. 
\end{theorem} 

We will now draw some consquences of this fundamental result. 

\begin{lemma} \label{mu through mu} 
The Sato microlocalization $Sh(X) \to Sh^{\RR^+}(T^*_M X)$ factors
through the (global sections of) a morphism of sheaves of categories 
\begin{equation}\label{micro factored through micro unconditional} 
\mu sh|_{T^*_M X} \to Sh^{\RR^+}
\end{equation}  
Here the right hand side is the (sheaf of categories of) conic sheaves on $T^*_M X$.

Moreover, for $\Lambda \subset T^*X$ a conic closed subset, this map restricts to
\begin{equation}\label{micro factored through micro} 
\mu sh_\Lambda|_{T^*_M X} \to Sh_{C_{T^*_M X}(\Lambda)}^{\RR^+}
\end{equation} 
where $C_{T^*_M X}(\Lambda)$ is the specialization of $\Lambda$ to the normal bundle
of $T^*M_X$, with that normal bundle then identified with $T^* (T^*_M X)$. 
\end{lemma} 
\begin{proof}
Because the target is a sheaf of categories, it is enough to construct a map from $\mu sh^{pre}$. 
To do this we should show that for any $\Omega \subset T^* X$, the microlocalization
$\mu_M$ induces a functor $Sh(X) / Sh_{T^*X \setminus \Omega}(X)  \to Sh(\Omega \cap T^*_M X)$.
In other words, we should show that if $F$ has no microsupport in $\Omega$, then 
$\mu_M(F)$ has no support in $\Omega \cap T^*_M X$.  This follows from Theorem 
\ref{microsupport of microlocalization}, which gives the microsupport of $\mu_M(F)$, so in particular
the support.  The final statement characterizing the behavior of microsupports is a direct translation
of Theorem \ref{microsupport of microlocalization}. 
\end{proof}

\begin{remark} 
Under the identification of $\mu hom$ with the Hom in $\mu sh$, the above functor acts
on Hom sheaves as the natural map

\begin{equation} \label{muhom to hom of mu} 
\mu hom (F, G)|_{T^*_M X} \to \mathcal{H} om(\mu_M F, \mu_M G)
\end{equation}
\end{remark} 

We are interested conditions on $\Lambda$ which ensure that (\ref{micro factored through micro}) is 
an equivalence, and in particular that the map (\ref{muhom to hom of mu}) is an isomorphism.  
It is possible to give counterexamples showing some condition is necessary for (\ref{muhom to hom of mu})
to be an isomorphism \cite{GS-letter}
(so in particular, the map (\ref{micro factored through micro unconditional}) is not an isomorphism).
One known sufficient condition is $\Lambda = T^*_M X$ (see e.g. \cite[Lem. 10.2.2, Prop. 10.2.4]{Gui2}, 
and \cite[Prop. 6.6.1, Lem. 7.5.2]{KS}, and more generally around \cite[Sec. 7.5]{KS}, for some
precursors).  This fact  is fundamental to the analysis 
of $\mu sh$ along a smooth conic Lagrangian.  

Here we note that, more generally, it suffices if $\Lambda$ is in an appropriate
sense {\em already} conic along $T^*_M X$, so the specialization to the normal cone is an innocent
operation. 

\begin{lemma} \label{mu is mu} 
Consider $\Lambda \subset T^*X$.  Suppose that, in the neighborhood of 
a point $\xi \in T^*_M X$, the set $\Lambda$ is contained in the union of conormals (in $T^*X$) 
to strata in a subanalytic (or otherwise o-minimal) stratification of $M$.  Then the map (\ref{micro factored through micro})
$$\mu sh_\Lambda|_{T^*_M X} \to Sh_{C_{T^*_M X}(\Lambda)}^{\RR^+}$$
is an isomorphism at $\xi$.  
\end{lemma}
\begin{proof}
For morphisms of sheaves valued in the category of categories, whether a map is an isomorphism
may be checked on stalks.  Thus the 
question reduces to a question about the $D(X; p)$ for $p \in T^*_M X$, to which the theory 
developed in \cite[Chap. 6, 7]{KS} applies.

We may replace $X$ by a tubular neighborhood of $M$, and fix a metric so as to identify this 
neighborhood with the normal bundle to $M$.  That is, we replace $X = T_M X$.   
Note that on sheaves conic with respect
to the scaling on $T_M X$, the given functor is just the Fourier-Sato transform, which is an equivalence.   

Now the point is just that the hypothesis of the theorem ensures that the sheaves in question
are microlocally conic, i.e. isomorphic to a conic on $T_M X$ sheaf in a neighborhood of 
$T^*_M X$.  That is because the conormals to strata in $M$ remain constant under the deformation
to the normal cone to $T^*_M X$.   Thus the Fourier transform is (microlocally) an isomorphism on these
sheaves.   
\end{proof}

\begin{remark}
Evidently the criterion of Lemma \ref{mu is mu} holds at every point away from 
the zero section in the FLTZ skeleton, when $M$ is taken
any sub-torus.  More generally, it is automatic away from the zero section 
for sheaves constructible with respect to a piecewise
linear stratification, where $M$ is amongst the strata.  
\end{remark}

\begin{cor} \label{cor:CCCmicro2}  
Let $\Sigma \subset M^\vee_\RR$ be a fan and $\LL_\Sigma$ the corresponding skeleton inside $T^* \TT^\vee_\RR$.

Let $\sigma \in \Sigma$ be a cone, and let $\pi: {\sigma^\perp} \times (-\sigma^\circ) \to \sigma^\perp$ be the projection.  Then the functor (from Lem. \ref{lem:CCCmicro}) 

\begin{eqnarray*}
m_\sigma: Sh_{\LL_\Sigma}(\TT^{\vee}) & \to & Sh_{\LL_{\Sigma / \sigma} }(\sigma^\perp) \\
F & \mapsto & \pi_* ((\mu_{\sigma^\perp} F)|_{\sigma^\perp \times (-\sigma^\circ)}).
\end{eqnarray*}

factors canonically through an isomorphism

$$\mu sh_{\LL_\Sigma}(\sigma^\perp \times (-\sigma^\circ)) \to Sh_{\LL_{\Sigma/\sigma}}(\sigma^\perp)$$

\end{cor}
\begin{proof}
As we have remarked, sheaves constructible with respect to a piecewise linear stratification 
necessarily satisfy the hypothesis of Lemma \ref{mu is mu}. 
\end{proof}

\subsection{At infinity}
We are now ready to pass to the boundary on both sides of Bondal's correspondence. On the B-side, this means passing 
from the toric variety $\bT_{\Sigma}$ to the union of its toric boundary divisors, and on the A-side, this means moving from 
the relative skeleton $\LL_{{\Sigma}}$ of the LG model $W: T^*\TT^\vee \to \CC$ to the complement of the zero section: 
$\LL_{\Sigma}^\circ:=\LL_{\Sigma}\setminus\TT^\vee.$

\begin{theorem}\label{thm:ccc-infty} For $\bT_\Sigma$ smooth, 
there is an equivalence of categories 
$$\Coh(\partial \bT_{\Sigma}) \cong \mu sh(\partial \LL_{\Sigma})^c.$$
\end{theorem}

\begin{proof}
To avoid worrying about whether various colimits exist, we will work with the cocomplete categories
IndCoh and $\mu sh$, and we will return
to the above statement at the end by passing to compact objects.  This is essentially only a matter of notation. 

According to Lemma \ref{lem:boundarypushout}, the toric boundary is a colimit of its component subvarieties.  By \cite[IV.4.A.1.2]{GR}, taking coherent
sheaves carries colimits to colimits (the result is there stated for affine schemes, from which the statement we require
follows by \'etale descent; note also that a colimit of underived schemes or stacks remains a colimit of the corresponding items
viewed as derived objects, since the inclusion of underived geometry into derived geometry is left adjoint to truncation of derived
structure).   Thus:
$
\IndCoh(\partial \bT_{\Sigma}) \cong \colim_{\sigma \in {\Sigma}} \IndCoh(\overline{O(\sigma)})$. 
By Zariski (or \'etale in the stack case) descent we may trade 
$\IndCoh(\overline{O(\sigma)}) \cong  
\Mod- B_\ocfan{\sigma}$.   (For a detailed explanation of this isomorphism  see \cite{Ku}.)

The coherent-constructible correspondence of \cite{B, FLTZ, Tr} and Kuwagaki's theorem \cite{Ku}, respectively, give the following two equivalences:
$$\colim_{\sigma \in {\Sigma}} \Mod-B_\ocfan{\sigma} \cong \colim_{\sigma \in {\Sigma}} \Mod-A_{\ocfan{\sigma}} \cong \colim_{\sigma \in {\Sigma}} 
 Sh_{\LL_{\ocfan{\sigma}}}(\TT(\sigma)^\vee).$$ 
Finally, by taking adjoints to the restriction morphisms we analyzed
in Lemma \ref{lem:CCCmicro} and Corollary \ref{cor:CCCmicro2}, we obtain the following identification: 
$$\colim_{\sigma \in {\Sigma}} 
 Sh_{\LL_{\ocfan{\sigma}}}(\TT(\sigma)^\vee) = \colim_{\sigma \in {\Sigma}}  \mu sh_{\LL_\Sigma}(\sigma^\perp \times (-\sigma^\circ))$$

On the right, the maps are co-restriction functors of wrapped microlocal sheaves, and this colimit is just the one associated to a cover 
of $\partial \LL_{\Sigma}$.  This completes the proof. 
\end{proof}

\begin{remark} \label{rem: blowdown}
The result holds without the smoothness hypothesis, as e.g. can be seen by taking some toric resolution, applying 
Theorem \ref{thm:ccc-infty}, and then matching the semiorthogonal decomposition of the category of the resolution 
on the $B$ side with stop removal on the $A$ side.  We content ourselves with the smooth case here because we anyway
have only in this case identified $\partial \LL_\Sigma$ as a Weinstein skeleton. 
\end{remark}

\section{A glimpse in the mirror of birational toric geometry}\label{sec:extra}

Since the works \cite{BO, K}, it has been understood that birational features of algebraic geometry often have natural
interpretations in the derived category of coherent sheaves.  
Mirror symmetry provides an illuminating perspective on these derived equivalences, which in algebraic geometry seem
to be among a discrete set of objects. Remarkably, on the mirror this discretization becomes unnatural, and one can
continuously interpolate between the mirrors of derived equivalent varieties. 
Many other features of birational geometry
(\emph{e.g.}, semi-orthogonal decompositions associated to blowups) also have beautiful new geometric interpretations in terms of mirror geometry. 
For discussions in the context of toric varieties, see \cite{FLTZ4, CKK, BDFKK}. 

Here is another result in this direction. 

\begin{cor}
Let $W:(\mathbb{C}^\times)^n\to \mathbb{C}$ a Laurent polynomial with Newton polytope $\Delta$ and 
$\Sigma_1,\Sigma_2$ a pair of fans obtained as star-shaped triangulations of $\Delta.$
Then there is a derived equivalence $\Coh(\bT_{\Sigma_1})\cong \Coh(\bT_{\Sigma_2}).$
	\label{cor:dereq}
\end{cor} 
\begin{proof}
Let $\LL_1, \LL_2$ be the corresponding \cite{FLTZ} skeleta.  We have shown that (a Liouville domain completing to)
the general fiber $\dtmir$ of $W$ is isotopic both to a domain with skeleton $\partial \LL_1$, and to a domain
with skeleton $\partial \LL_2$.  By \cite[Cor. 2.9]{GPS2}, we have an equivalence of the wrapped Fukaya categories
$Fuk(T^* \TT^\vee, \partial \L_1) \cong Fuk(T^* \TT^\vee, \partial \L_2)$.  By \cite{GPS3} we may trade this
for an equivalence of constructible sheaf categories, and by \cite{Ku} we may trade the latter for the asserted equivalence
of coherent sheaf categories. 
\end{proof}

What the above argument does not yet give is a formula for the above equivalence.  In fact, there are many such 
derived equivalences, corresponding to monodromies (as we vary the coefficients of $f$) around the discriminant 
locus. We will describe these in future work.

\subsubsection{Non-Fano mirror symmetry}
It was observed in \cite{AKO, A2} 
that mirror symmetry for toric varieties requires modification in the case of a non-Fano variety $\bT$: the naive interpretation
of the Hori-Vafa mirror Landau-Ginzburg model for a non-Fano variety contains $\Coh(\bT)$ as a full subcategory but can be strictly larger. 
One procedure to remedy this discrepancy is suggested in \cite{BDFKK}.  By contrast, \cite{Ku} holds for {\em all} toric varieties.  
Here we explain this discrepancy in an example; in future work, we plan to use the same ideas to establish
the conjectures of \cite{BDFKK}. 

In the body of this paper, we began with a polytope $\Delta^\vee$ with star-shaped triangulation, and let $\Sigma$ be the fan given by this star-shaped triangulation.
Any fan $\Sigma$ obtained in this way has the following property: let $v_1,\ldots,v_k$ be (stacky) primitives for the rays in $\Sigma$, and let $\Delta^\vee$ be the convex hull of the $v_i$. Then each $v_i$ is on the boundary of $\Delta^\vee.$
This excludes fans $\Sigma$ in which one of the primitives $v_i$ is too short to reach $\partial\Delta^\vee.$ In this case, the the mirror to $\partial\bT_\Sigma$ will not be a hypersurface with Newton polytope $\Delta$, but only a Liouville subdomain of such a hypersurface.
The simplest case of this is described in the following example.

\begin{example}
	Let $\Sigma_1$ be the fan with primitive rays $(-1,2), (1,2)$;
	$\Sigma_2$ the fan with primitives $(-1,2), (0,2), (1,2)$;
	$\Sigma'$ the fan with primitives $(-1,2), (0,1), (1,2)$;
	and $\Delta^\vee$ the polytope obtained as convex hull of the primitives for any of the three fans above. (These convex hulls obviously agree.)
	\begin{figure}[h]
		\centering
		\includegraphics[width=10cm]{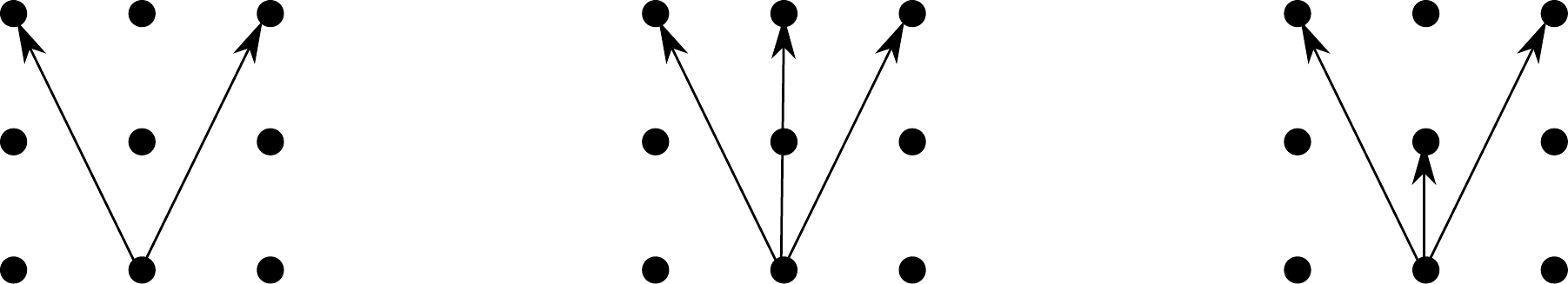}
		\caption{The fans $\Sigma_1,\Sigma_2,\Sigma'$.}
		\label{fig:nonfano-fans}
	\end{figure}
	Then each of $\Sigma_1$ and $\Sigma_2$ is obtained as a star-shaped triangulation of $\Delta^\vee$; 
	hence the results of this paper show that the boundaries $\partial \bT_{\Sigma_1}$ and $\partial \bT_{\Sigma_2}$ are both mirror to a generic hypersurface $H$ with Newton polytope $\Delta^\vee.$
	\begin{figure}[h]
		\centering
		\includegraphics[width=10cm]{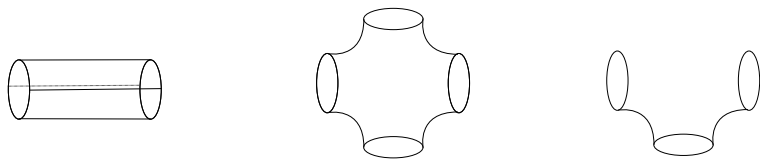}
		\caption{The FLTZ boundary skeleta $\Lambda_1,\Lambda_2,\Lambda'$ for the fans $\Sigma_1,\Sigma_2,\Sigma'$.}
		\label{fig:nonfano-skeleta}
	\end{figure}

	We obtain two different skeleta $\Lambda_1,\Lambda_2$ of the hypersurface $H$, corresponding to the respective triangulations $\Sigma_1$ and $\Sigma_2$, and
	we know that each of these is the boundary of a stacky FLTZ skeleton; by studying the fans $\Sigma_i$, we conclude that $\Lambda_1$ consists of two circles connected by four different intervals (since the two rays in $\Sigma_1$ share a non-unimodular simplex of area 4), and $\Lambda_2$ consists of four circles, cyclically connected by intervals (there being four circles since the middle ray, of length two, is double-counted by the stacky FLTZ procedure).
	Each of these is a skeleton for $H$, which is a quadruply-punctured genus-one curve.

	Let $\Lambda'$ be the boundary of the \cite{FLTZ} skeleton for $\Sigma'$. Then $\Lambda'$ is no longer a skeleton for the hypersurface $H$, as $\Lambda_1,\Lambda_2$ are. It resembles the skeleton $\Lambda_2$, except that the central ray, now of length one, is no longer double-counted. This means that $\Lambda'$ is obtained from $\Lambda_2$ by deleting one of the two double-counted circles along with its two connecting intervals. Hence $\Lambda'$ consists of three circles, connected in a row by a pair of intervals. It is the skeleton of a triply-punctured genus-one curve, a subdomain of the quadruply-punctured curve $H$.
\end{example}

\end{document}